\documentclass[11pt]{article}
\usepackage{abstract}
\usepackage{makecell}
\usepackage{array}
\usepackage{wrapfig}
\usepackage[noadjust]{cite}
\usepackage[utf8]{inputenc}
\usepackage{authblk}
\usepackage{amsmath}
\usepackage{blindtext}

\usepackage{color}

\newcommand{\RM}[1]{{#1}}

\usepackage{hyperref}
\usepackage{url}
\usepackage{geometry}
 \usepackage{booktabs}
\usepackage{amsfonts} 
\usepackage{amssymb}
\usepackage{float}
\usepackage{}
\usepackage[nottoc]{tocbibind}
\usepackage{graphicx}
\graphicspath{{./image/}}
\usepackage{mathtools}
\renewenvironment{abstract}
 {\par\noindent\textbf{\abstractname:}\ \ignorespaces}
 {\par\medskip}
\usepackage{multicol}
\font\myfont=cmr12 at 12pt
\title{Classification and stability analysis of polarising and depolarising travelling wave solutions for a model of collective cell migration}

\author{\large{ \textbf{Nizhum Rahman$^{1,2,\ast}$, Robert Marangell$^{3}$ and Dietmar Oelz$^{1}$} }}
\affil{\myfont $^{1}$ School of Mathematics and Physics, The University of Queensland, QLD, 4072, Australia. \vspace{-2ex}}
\affil{\myfont $^{2}$ Faculty of Science and Information Technology, Daffodil International University Dhaka,1207, Bangladesh. \vspace{-2ex}}
\affil{\myfont $^{3}$ School of Mathematics and Statistics, University of Sydney, NSW, 2006, Australia. \vspace{-2ex}}
\begin{document}
\maketitle
\begin{abstract}

We study travelling wave solutions of a 1D continuum model for collective cell migration in which cells are characterised by position and polarity. Four different types of travelling wave solutions are identified which represent polarisation and depolarisation waves resulting from either colliding or departing cell sheets as observed in model wound experiments. 
We study the linear stability of the travelling wave solutions numerically and using spectral theory. This involves the computation of the Evans function most of which we are able to carry out explicitly, with one final step left to numerical simulation.
\newline 

\end{abstract}

\section{Introduction}
Epithelial cells line the surfaces of our body. They are tightly packed and found in various organ systems.
They serve many roles such as secretion, absorption, sensation, protection and transport. In the context of various physiological and pathological processes including tissue repair, cancer and wound healing \cite{a1,a2} they undergo collective cell migration which is an important feature characterising the development and life-cycle of multi-cellular organisms.
The analysis of mathematical models for cell migration is a powerful tool that allows to identify characteristic features of the collective dynamics.

Cells can move as a group with shared responsibilities according to recent studies on migratory epithelial tissues \cite {a3,a4,a5,a6,a7}. During collective cell migration epithelial cells maintain stable cell-cell junctions \cite{a8,a9,a10,a11}. 

Directed cell migration is closely linked to the polarity of the cells \cite{a12}, which refers to the spatial differences between cells in shape, structure and function. 

How cell polarity is influenced by intercellular coupling is a fundamental question. The best way of understanding its underlying biochemical signalling mechanism is the planar cell polarity pathway coupling bistable intracellular states among adjacent cells \cite{a13,a14,a15,a16,Oelz2019}. 

In cells, polarity propagates as a travelling wave from cell to cell \cite{a17,a18,Oelz2019}. The corresponding travelling wave solution of a 1D model for the collective migration of epithelial cells has been computed explicitly in \cite{Oelz2019}.
Preliminary numerical simulations in \cite{Oelz2019} have shown that the same model exhibits other travelling wave solutions such as a depolarisation wave entering the sheet of departing cells. 

In the present article we identify various other types of travelling wave solutions of this model. We analyse their stability with the help of spectral theory and using tools from dynamical systems theory \cite{SANDSTEDE2002,a20,a21} such as the Evans function which provides information about the point spectrum of the linearised model. Note that we are able to compute most of the information to setup the Evans function explicitly with only one last component being evaluated numerically. 

This article is structured as follows. Section \ref{sec_prel} states the governing equations of the 1d model for collective cell migration. Section~\ref{sec_classification} lists various types of polarising and depolarising travelling wave solutions which can be formulated building on the solution found in \cite{Oelz2019}. In section~\ref{sec_numerics} we provide numerical evidence for their stability and in section~\ref{sec_stabanalysis} we perform the linear stability analysis. Finally we wrap up and discuss these results in section~\ref{sec_concl}.

\section{Preliminary results}
\label{sec_prel}

\begin{figure}[H]
	\centering
	\includegraphics[width=0.9\linewidth]{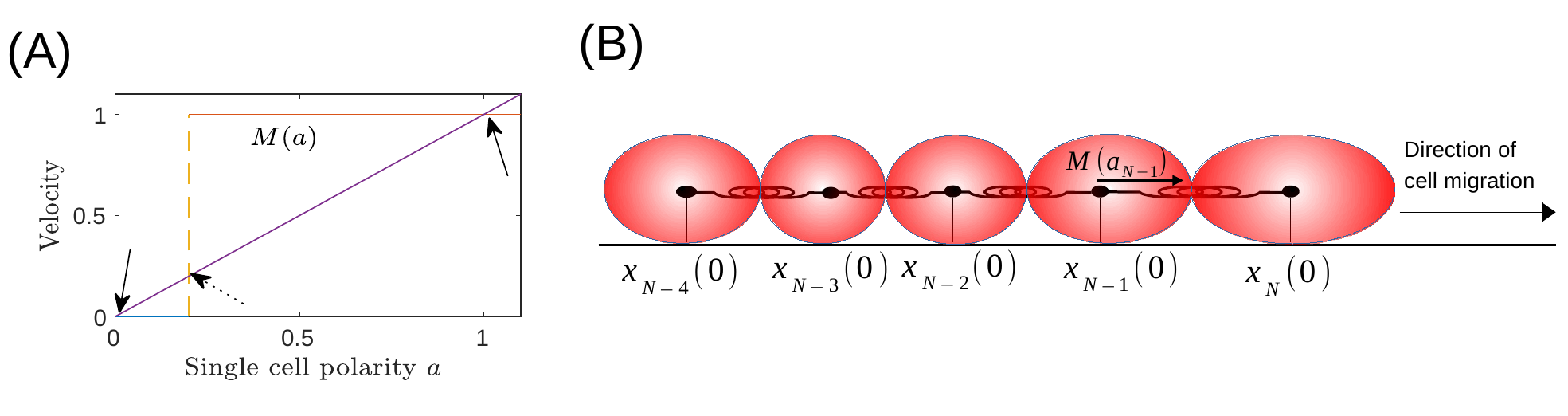}
	\caption{A: Visualisation of the steady states (arrows) of the single cell model defined by $M(a)=a$. B: Sketch of 1d model for collective cell migration.}
	\label{fig_singlecell}
\end{figure}

The mathematical model for collective migration of epithelial cells introduced in \cite{Oelz2019} characterises each cell by its position along the real axis and by a real-valued quantity called polarity $a \in \mathbb{R}$. It represents the asymmetry of the cell and it is assumed to be linked to active migratory velocity through a non-linear function $v_i=M(a_i)$.
This function may be thought of as the following step function (or a smoothened version of it as in Fig.~\ref{fig_singlecell}A)
\begin{equation}
\label{equ_M}
M(a)=\begin{cases}
0 & a \leq \alpha \; , \\
1 & a > \alpha  \; ,
\end{cases}
\end{equation}
which models polarity-dependent motility with threshold polarity $\alpha>0$.

For a single cell the following differential equation for the polarity is defined which models auto-depolarisation of the cell and adaption to actual motion. It is given (omitting physical constants) by $\dot a=-a + v=-a+M(a)$ which - taken as a dynamical system - features two stable steady states, a non-polarised one at $a=0$ and a polarised one at $a=1$ (see Fig.~\ref{fig_singlecell}).

The 1d model introduced in \cite{Oelz2019} treats the epithelial cell sheet as a 1d chain of such cells connected by linearly elastic springs.
As a consequence each cell's velocity is determined by the spring forces emanating from neighbouring cells and by the active migratory force. The governing equations (after non-dimensionalisation) are 
\begin{equation}
\label{equ_particle}
\left\{
\begin{aligned}
&\dot a_i = - a_i + \dot x_i \; , \\
&\dot x_i =M(a_i) + \kappa (x_{i_1}-2 x_i + x_{i-1}) \; ,
\end{aligned}
\right.
\end{equation}
where $\kappa$ is a phenomenological parameter representing cellular contractility.

The continuum model associated to this particle model has been derived in \cite{Oelz2019}. Its  governing equations are the continuity equation for the cell density $\rho=\rho(x,t)$ coupled to the momentum equation and a separate equation for the cell polarity $a(x,t)$ which correspond to the system \eqref{equ_particle}. After non-dimensionalisation, these equations are given by
\begin{equation}
\label{equ_main}
     \begin{cases}
    \frac {\partial \rho}{\partial t}+\frac{\partial}{\partial x}(\rho v)=0,\\
    \frac{\partial a}{\partial t}+v\frac{\partial a}{\partial x}=-a+v,\\
    v= M(a)-\kappa\frac{1}{\rho^3} \frac{\partial\rho}{\partial x},
    \end{cases}
\end{equation}
which involves the velocity field as $v=v(x,t)$ and the cell polarity $a=a(x,t)$. 

In \cite{Oelz2019} the system \eqref{equ_main} is reformulated using the travelling wave ansatz  $\rho(x,t)=R(z),a(x,t)=A(z)$, where $z = x - st$ parametrises the wave profiles and $s$ is the wave speed. Then the reduced system of equations is coupled to boundary conditions which correspond to one of the stable fixed points of the single particle model \eqref{equ_particle} prescribed at $z=-\infty$, namely $( R_{-\infty},A_{-\infty})=(1,0)$, where $R_{-\infty}$ is the cell density at rest normalised to 1. This yields
\begin{equation}\label{equ_travwave}
     \begin{cases}
    R'=\frac{R^2}{\kappa}((M(A)-s)\, R +s),\\
    A'=1+(\frac{A-s}{s}) \, R. 
    \end{cases}
\end{equation}
Note that the second stationary point of \eqref{equ_travwave} is $(R_\infty,A_\infty)=(\frac{s}{s-1},1)$.

\begin{figure}[H]
	\centering
    \includegraphics[width=.99\textwidth]{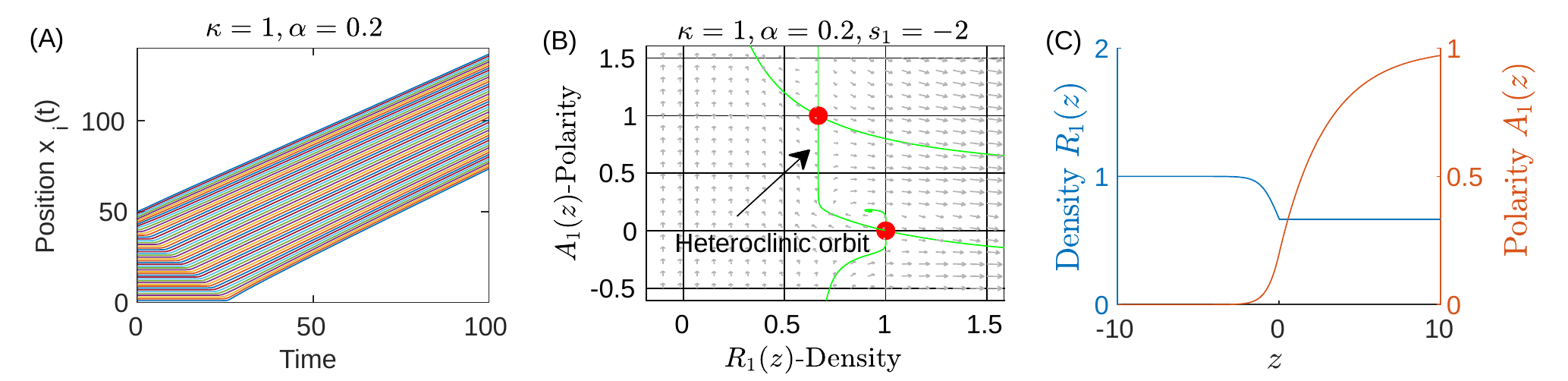}
    \caption{(A) shows a simulation of \eqref{equ_particle} corresponding to the  polarisation caused by a departing cell sheet. (B) shows the corresponding phase plane diagram of \eqref{equ_travwave} and (C) shows the travelling wave profiles of density R(z) and polarity A(z).}
    \label{fig_wavep1}
\end{figure}

In \cite{Oelz2019} the authors identify an explicit solution to the travelling wave problem for the specific active velocity function given in \eqref{equ_M}.
The travelling wave solution they single out models ongoing polarisation of cells initiated by a departing cell sheet (Fig.~\ref{fig_wavep1}). Its travelling wave speed and profiles $R$ and $A$ (all with subscript $1$) are 
\begin{enumerate}
	\item [S1:] Polarisation wave triggered by departing cell (Fig~\ref{fig_wavep1}) sheet\cite{Oelz2019}: $s_1=-\sqrt{\kappa \left(\frac{1}{\alpha}-1 \right)}$ with travelling wave profiles
	\begin{equation*}
	A_1(z)=\begin{cases}
	s_1\alpha\left(1-\frac{1}{g^{-1}(M-\frac{zs_{1}}{\kappa})}\right) \; ,& z <\mbox {0}\\
	1+(\alpha-1)e^{\frac{z}{s_1-1}}\;, & z \geq\mbox{0}
	\end{cases}
	\quad \text{and} \quad
	R_1(z)=\begin{cases}g^{-1}\left(M-\frac{zs_1}{\kappa}\right)\;,& z <\mbox{0}\\
	\frac{s_1}{s_1-1}\;, & z \geq\mbox{0} \; 
	\end{cases}
	\end{equation*}
	where $g(y)=\frac{1}{y}+\log\left(\frac{1}{y}-1\right)$ and $M= g\left(\frac{s_1}{s_1-1}\right)$.
\end{enumerate}

\section{Travelling wave solutions}
\label{sec_classification}
In addition to the travelling wave solution S1 found in \cite{Oelz2019} we identify other travelling wave solutions corresponding to polarization and depolarization waves. The specific travelling wave speeds and profiles may be obtained through computations which are analogous to those performed in \cite{Oelz2019}. Alternatively, one may also derive them directly through two different transformations. These specific transformations of dependent and independent variables allow us to recast the travelling wave solution S1 into other travelling wave solutions.

\begin{figure}[H]
    
    \centering
    \includegraphics[width=0.7\textwidth]{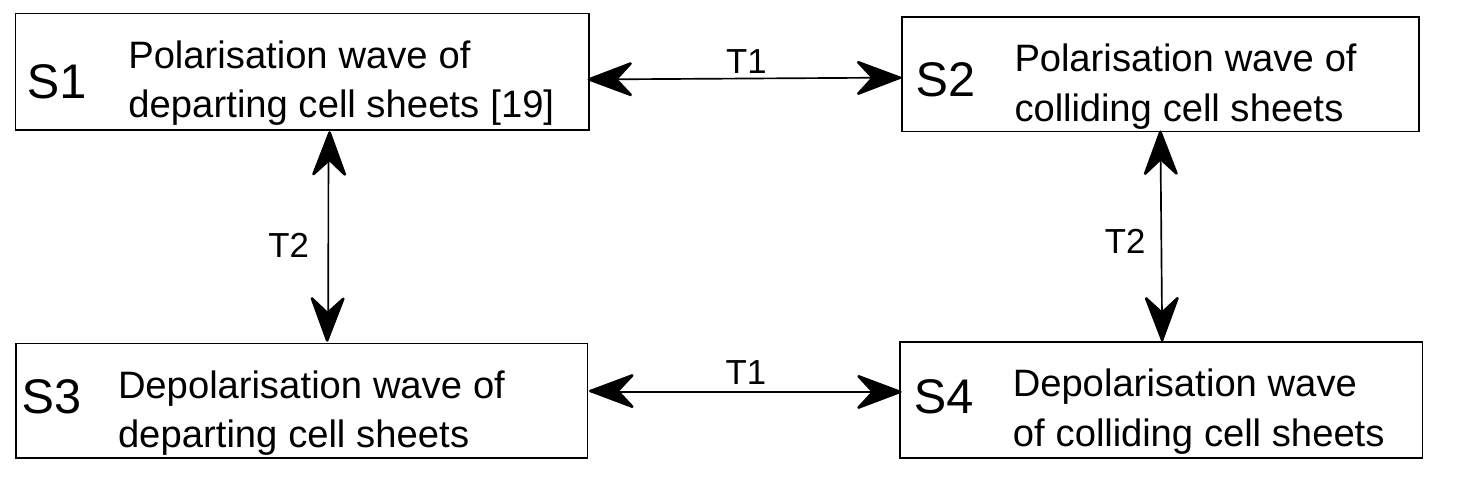}
    \caption{Scheme of travelling wave solutions (polarising/depolarising, departing/colliding) obtain by applied the transfromations $T_1$ and $T_2$.}
    \label{Fig_transforms}
\end{figure}

We consider two different transformations:
\begin{enumerate}
     \item 
 One which transforms a solution of system \eqref{equ_travwave} into another solution of the same system,
\begin{equation*}
    T_1:(R(z),A(z), s)\longrightarrow(\Bar{R}(\Bar{z}),\Bar{A}(\Bar{z}), \Bar{s}),
\end{equation*}
where $\Bar{A}=A$, $\frac{1}{R}+\frac{1}{\Bar{R}}=2$, $\Bar{z}=
\int_{0}^{z} \left(1-2 \, R(\tilde{z}) \right) d\tilde{z}$ and $\bar s=-s$,\newline
\item 
and a transformation which transforms a solution to \eqref{equ_main} to another solution of \eqref{equ_main},
\begin{equation*}
    T_2:(\rho(x,t),a(x,t))\longrightarrow({\Bar{\rho}}({\Bar{x}},{\Bar{t}}),{\Bar{a}}({\Bar{x}},{\Bar{t}})),
    \end{equation*}
    where ${\Bar{t}}=t $, ${\Bar{x}}=t-x$, ${\Bar{v}}=1-v$, ${\Bar{a}}=1-a$, $\Bar{\alpha}=1-\alpha$ and $\Bar{\rho}=\rho$.

    Specifically for travelling wave solutions the transformation $T_2$ is given by
    \begin{equation*}
    	\tilde T_2:(R(z),A(z), s)\longrightarrow(\Bar{R}(\Bar{z}),\Bar{A}(\Bar{z}), \Bar{s}) \; , \quad 
    	\text{where} \quad \Bar{A}=1-A \; ,  \quad \Bar{R}=R \quad \text{and}\quad \Bar{s}=1-s \; .
    \end{equation*}
 \end{enumerate}
Applying these transformations, we obtain three additional travelling wave solutions in addition to S1 as illustrated in Fig.~\ref{Fig_transforms}.

\begin{figure}[H]
	\centering
    \includegraphics[width=0.99\textwidth]{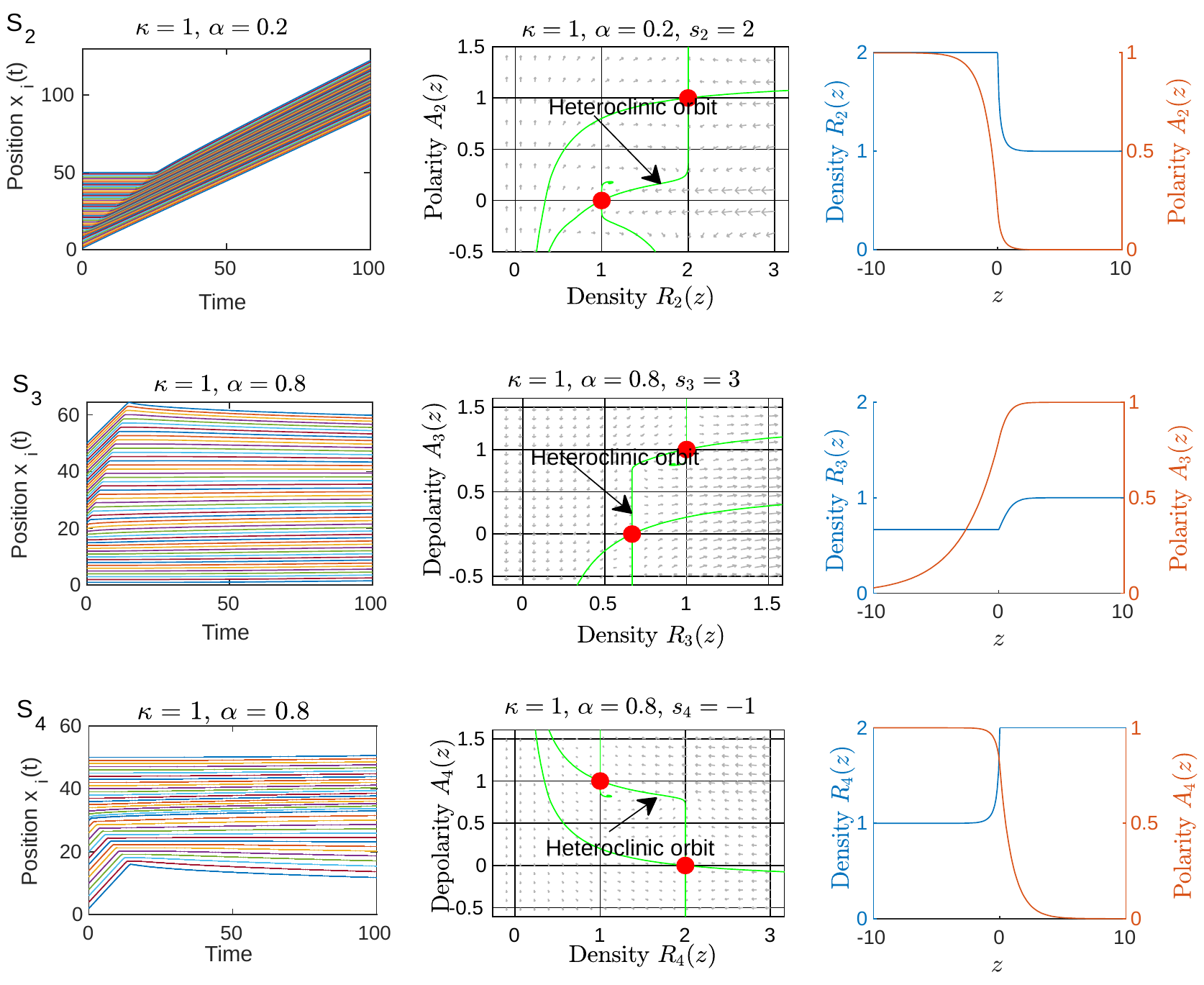}
    \caption{S2: Polarisation wave caused by a colliding cell sheet. S3: Depolarisation wave due to a departing cell sheet. S4: Depolarisation wave caused by a colliding cell sheet. The second and third columns of S1, S2 and S3 show the corresponding phase plane diagrams of solutions to \eqref{equ_travwave} and the travelling wave profiles respectively.}
    \label{fig:wavep2_3_4}
\end{figure}

Their specific speeds and profiles are:
\begin{enumerate}
    \item [S2:] By applying the transformation $T_1$ to the travelling wave solution (3), we obtain the polarisation wave caused by a colliding cell sheet. Its speed is  given by $s_2=\sqrt{\kappa\left(\frac{1}{\alpha}-1\right)}$ and the wave profiles $A_2$ and $R_2$ (third column of Fig.~\ref{fig:wavep2_3_4}~S2) are given by
\begin{equation*}
 A_2(z)=\begin{cases}
s_2\alpha \left(1-\frac{1}{h^{-1}(N-\frac{zs_2}{\kappa})}\right)\;,& z \geq\mbox {0}\\
1+(\alpha-1)e^{\frac{z}{s_2-1}}\;, & z <\mbox{0}
\end{cases}
\quad \text{and} \quad 
R_2(z)=\begin{cases}h^{-1}\left( N-\frac{zs_2}{\kappa}\right)\;, & z \geq\mbox{0}\\
\frac{s_2}{s_2-1}\;, & z <\mbox{0}
\end{cases}
\end{equation*}
where $h(y)=\frac{1}{y}+\log(1-\frac{1}{y})$ 
and $N=h\left(\frac{s_2}{s_2-1}\right)$. Note that the travelling wave solution S2 can only be realised if the model parameters $\kappa$ and $\alpha$ are such that its wave speed $s_2 > 1$. If that is not the case the mathematical solution is not physical and violates the impenetrability of single cells (see supplementary section~\ref{sec_unphysical}).
\item [S3:]
After applying $T_2$ to the solution S1, we find the depolarization wave due to a departing cell sheet. The profiles $A_3$ and $R_3$ (third column of Fig.~\ref{fig:wavep2_3_4}~S3) are given by
\begin{footnotesize}
\begin{equation*}
A_3(z)=\begin{cases}
1+\left(1-s_3\right)\left(1-\alpha\right)\left(\frac{1}{g^{-1}\left(M+\frac{(1-s_3)z}{\kappa}\right)}-1 \right),& z \geq\mbox {0}\\
\alpha e^{\frac{z}{s_3}}, & z <\mbox{0}
\end{cases}
\; \text{and}\; 
R_3(z)=\begin{cases}g^{-1}\left(M+\frac{(1-s_3)z}{\kappa}\right),& z \geq\mbox{0}\\
\frac{s_3-1}{s_3}, & z <\mbox{0}
\end{cases}
\end{equation*}
\end{footnotesize}
where $M= g\left(\frac{s_3-1}{s_3}\right)$ and the wave speed is given by $s_3=1+\sqrt{\kappa\left(\frac{1}{1-\alpha}-1\right)}$.
\item[S4:]
After applying $T_2$ to the solution S2, we find the depolarization wave caused by a colliding cell sheet. The profiles $A_4$ and $R_4$ (third column of Fig.~\ref{fig:wavep2_3_4}~ S4) are given by
\begin{footnotesize}
\begin{equation*}
A_4(z)=\begin{cases}
1+\left(1-s_4\right)\left(1-\alpha\right)\left(\frac{1}{h^{-1}\left(N+\frac{(1-s_4)z}{\kappa}\right)}-1\right),& z <\mbox {0}\\
\alpha e^{\frac{z}{s_4}}, & z \geq\mbox{0}
\end{cases}
\; \text{and} \;
R_4(z)=\begin{cases}h^{-1}\left(N+\frac{(1-s_4)z}{\kappa}\right),& z <\mbox{0}\\
\frac{s_4-1}{s_4}, & z \geq\mbox{0}
\end{cases}
\end{equation*}
\end{footnotesize}
where $N= h\left(\frac{s_4-1}{s_4}\right)$ and the wave speed is given by $s_4=1-\sqrt{\kappa\left(\frac{1}{1-\alpha}-1\right)}$. Note that similar to S2 the travelling wave solution S4 violates the impenetrability of cells if the model parameters $\kappa$ and $\alpha$ are such that the wave speed satisfies $s_4 \geq 0$. (see supplementary section~\ref{sec_unphysical}).
\end{enumerate}
In the rest of the paper we will be concerned with the stability of the travelling wave solutions $S1$-$S4$. Note that the transformation $T_2$ acts on solutions of the continuum model \eqref{equ_main} and in a vicinity of the travelling wave solution it is a smooth map. Therefore any perturbation of one the travelling wave solutions S1 and S2, respectively, translates into a perturbation of S3 and S4, respectively - and vice versa. Therefore we will only investigate the spectral stability of the travelling wave solution S1, which will immediately imply the spectral stability of S3. 

The same argument will also apply to S2 and S4, respectively. Yet, since the transformation T1 only applies to solutions of \eqref{equ_travwave}, but not of \eqref{equ_main} we have to investigate the stability of S2 separately along the same lines as for S1. For this reason this will be added in the supplementary section~\ref{sec_stabanalysisS2}.

\section{Numerical Simulation}
\label{sec_numerics}

To explore the stability of the travelling wave solution S1 we compute numerical solutions of \eqref{equ_main} numerically using a Lax-Friedrichs scheme \cite{lax}. Typically we use very fine spatial grids to minimise the approximation error. 

\begin{figure}[H]
\centering
    \includegraphics[width=1\textwidth]{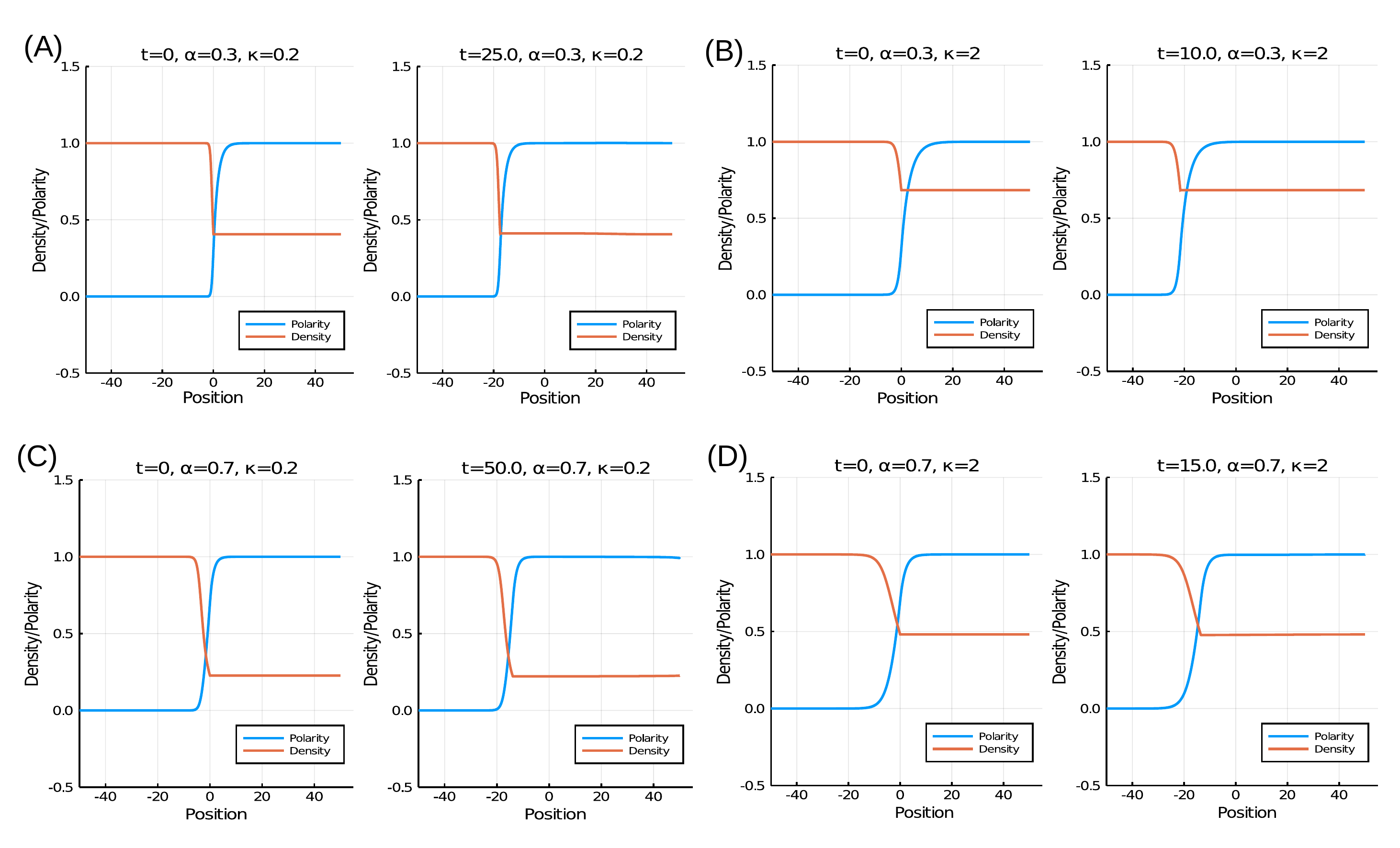}
    \caption{Travelling wave solution S1 for different parameter values $\alpha$ and $\kappa$.}
    \label{numerical_travellingwaves1}
\end{figure}

We are interested in the question whether the choice of the threshold polarity $0< \alpha < 1$ which models the sensitivity of cells to polarisation affects the stability of the travelling wave solutions.
Running the simulations starting with the initial condition given by the travelling wave solution S1 shows that for both, small values of $\alpha$ and large values of $\alpha$, the travelling wave solution is stable (Fig~\ref{numerical_travellingwaves1}(A,B) and Fig~\ref{numerical_travellingwaves1}(C,D)).

We take this as an indication that that the travelling wave S1 (and therefore S3) is stable for all parameter values. Further below we will investigate this question using spectral analysis of the linearised operator.

\begin{figure}[H]
	\centering
	\includegraphics[width=\textwidth]{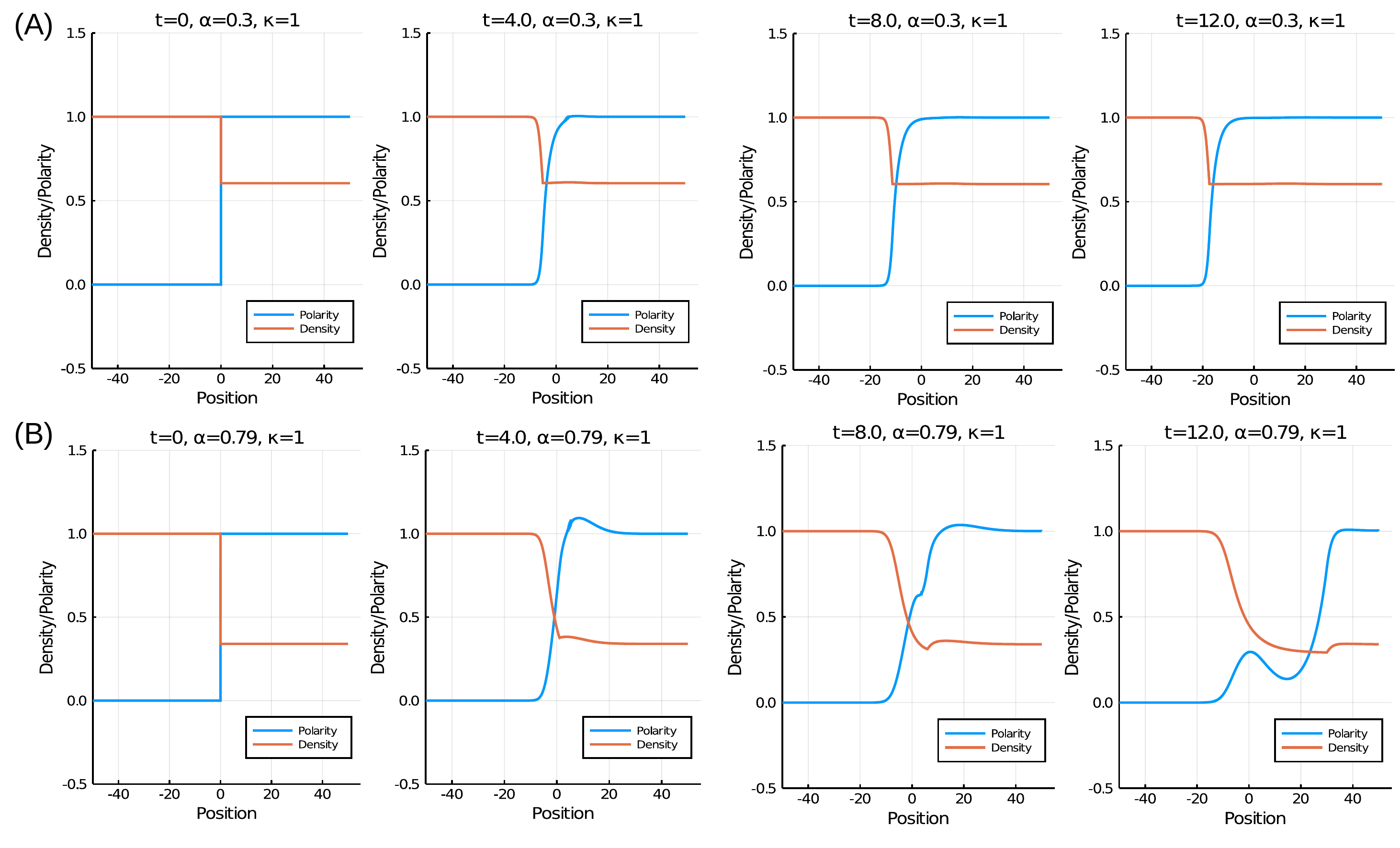}
	\caption{(A): Travelling wave for small $\alpha$. (B): Depolarization travelling wave due to large $\alpha$.}
	\label{small_alpha}
\end{figure}

\begin{figure}[H]
	\centering
	\includegraphics[width=.6\textwidth]{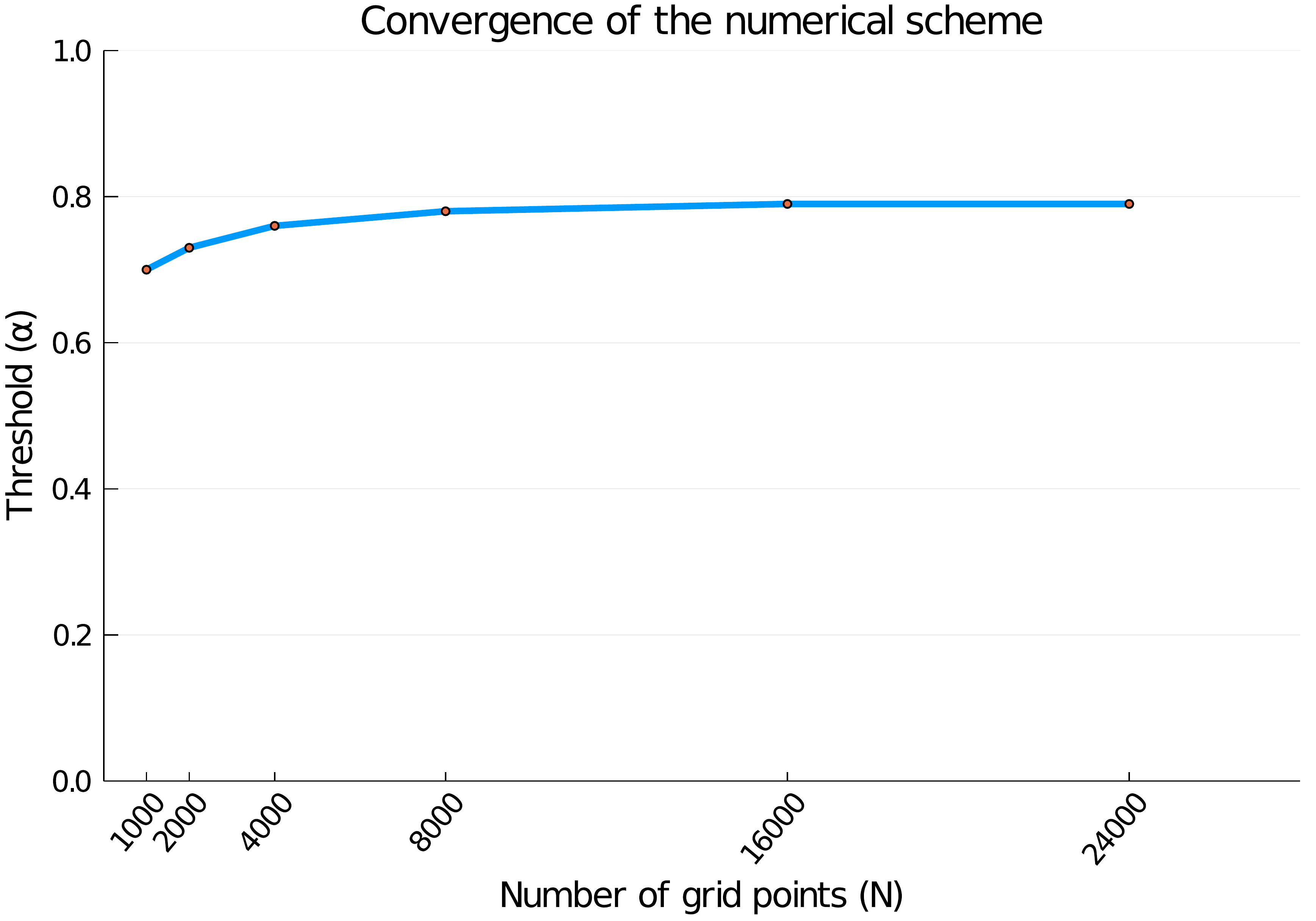}
	\caption{The threshold value $\bar \alpha$ such that the solution converges to the polarisation wave if $\alpha<\bar \alpha$, and to the depolarisation wave otherwise.}
	\label{Threshold}
\end{figure}

We also perform a numerical experiment in order to explore whether the size of the attractive regions of the polarisation wave S1 and the respective depolarisation wave S3 may depend on the sensitivity $\alpha$. 
To this end we simulate \eqref{equ_main} starting with a given initial condition that is distinct from the two travelling wave profiles, namely a step function centred at $x=0$. The simulations for small $\alpha$ (Fig~\ref{small_alpha}(A)) illustrate the convergence of the solution towards the polarisation wave. The same simulation for a large value of $\alpha$, however,  (Fig~\ref{small_alpha}(B)) shows that the solution readjusts and exhibits a depolarisation wave in the opposite direction. Taking into account potential approximation errors due to the numerical discretisation, we find that the threshold value $\bar \alpha$ (i.e. convergence to polarisation wave if $\alpha<\bar \alpha$, otherwise convergence to depolarisation wave) converges to about $\bar \alpha =0.79$ as the spatial grid is gradually tuned finer (Fig~\ref{Threshold}).

\section{Stability Analysis}
\label{sec_stabanalysis}
We now investigate the spectral stability of these travelling wave solutions. To this end we introduce the perturbations of density $\delta \rho$, polarity $\delta a$ and velocity $\delta v$ and linearise the system \eqref{equ_main} in the moving coordinate frame $(t, z=x-s t)$ at the travelling wave profile S1, i.e. $s=s_1$ and in what follows we write $R=R_1$, $A=A_1$ and $V$ for the associated velocity field. The linearised system of equations is given by
\begin{equation}
\label{equ_linearised}
\begin{aligned}
\begin{pmatrix}
\partial_t \delta \rho\\
\partial_t \delta a
\end{pmatrix} & = \mathcal{L} \begin{pmatrix}
\delta \rho\\
\delta a
\end{pmatrix} = \begin{pmatrix}
s \delta\rho'-(V \delta \rho+R \delta v)'\\
s\delta a'-\delta v( A'-1)- V \delta a'-\delta a 
\end{pmatrix}\;, \\
&\text{where} \quad \delta v =M'(A)\delta a-\kappa\left(\frac{1}{R^3}\delta\rho \right)' \; ,
\end{aligned}
\end{equation}
and where we use the notation $'=\partial_z$.
To obtain the associated eigenvalue problem $(\mathcal{L}-\lambda) (\delta \rho, \delta a)^T=0$, consider the solution of the above system \eqref{equ_linearised} to be $\delta \rho(z,t)=e^{\lambda t} \delta \rho(z)$ and $\delta a(z,t)=e^{\lambda t} \delta a(z).$ This converts the eigenvalue problem into a system of first order ODEs,
\begin{equation}
\label{final_A}
 \begin{pmatrix}
 \delta\rho'\\ \delta v'\\   \delta a'
 \end{pmatrix} = A(z,\lambda)\begin{pmatrix}
 \delta\rho\\     \delta v\\           \delta a     \end{pmatrix}\;,
 \quad \text{where}\quad
A(z,\lambda)=\begin{pmatrix}                         \frac{3\,R'}{R}&\frac{-R^3}{\kappa}&\frac{M'\,R^3}{\kappa}\\            \frac{2\,R's-R^2\,\lambda}{R^3}&\frac{-R'}{R}-\frac{R\,s}{\kappa}&\frac{M'\,R\,s}{\kappa}\\
0 &\frac{ (A'-1)\,R}{s}& \frac{R\,(\lambda+1)}{s}
\end{pmatrix} \; ,
\end{equation}

where we used that the velocity field (third equation in \eqref{equ_main}) for the travelling wave solution S1 satisfies 
\begin{equation}
\label{equ_V}
V=s(1-1/R) \quad \text{and} \quad V'=s\,R'/R^2 \; . 
\end{equation}

We then denote the corresponding linear operator  by $\mathcal{T}$, where $\mathcal{T}(y)(z)=(\frac{\partial}{\partial z}-A(z,\lambda))y$ and $y=(\delta\rho, \delta v, \delta a)^{T}$.
\par                          
It is our goal to find the spectrum of the linearized operator  $\mathcal{L} $ which is an operator on $H^2\times H^1$ mapping into $L^2\times L^2$ \cite{SANDSTEDE2002,a24}. If $(\mathcal{L}-\lambda)^{-1}$ is unbounded or does not exist for $\lambda\in\mathbb{C}$, then $\lambda$ is in the spectrum $\sigma(\mathcal{L})$. The Fredholm index of $\mathcal{L}$ is $ind(\mathcal{L}) =dim[ker(\mathcal{L})]-codim[R(\mathcal{L})]$, where $R(\mathcal{L})$ and $ker(\mathcal{L})$ denote the range and kernel  of $\mathcal{L}$ respectively \cite{a24}.  The spectrum of a Fredholm operator $\mathcal{L}$  is composed of two disjoint sets, namely point spectrum and essential spectrum. The point spectrum will consist of values $\lambda\in\sigma(\mathcal{L})$  such that $(\mathcal{L}-\lambda)$ is a Fredholm operator of index zero. The essential spectrum is the complement of the point spectrum in $\sigma(\mathcal{L})$.

\begin{figure}
	\centering
	\includegraphics[width=0.45 \textwidth]{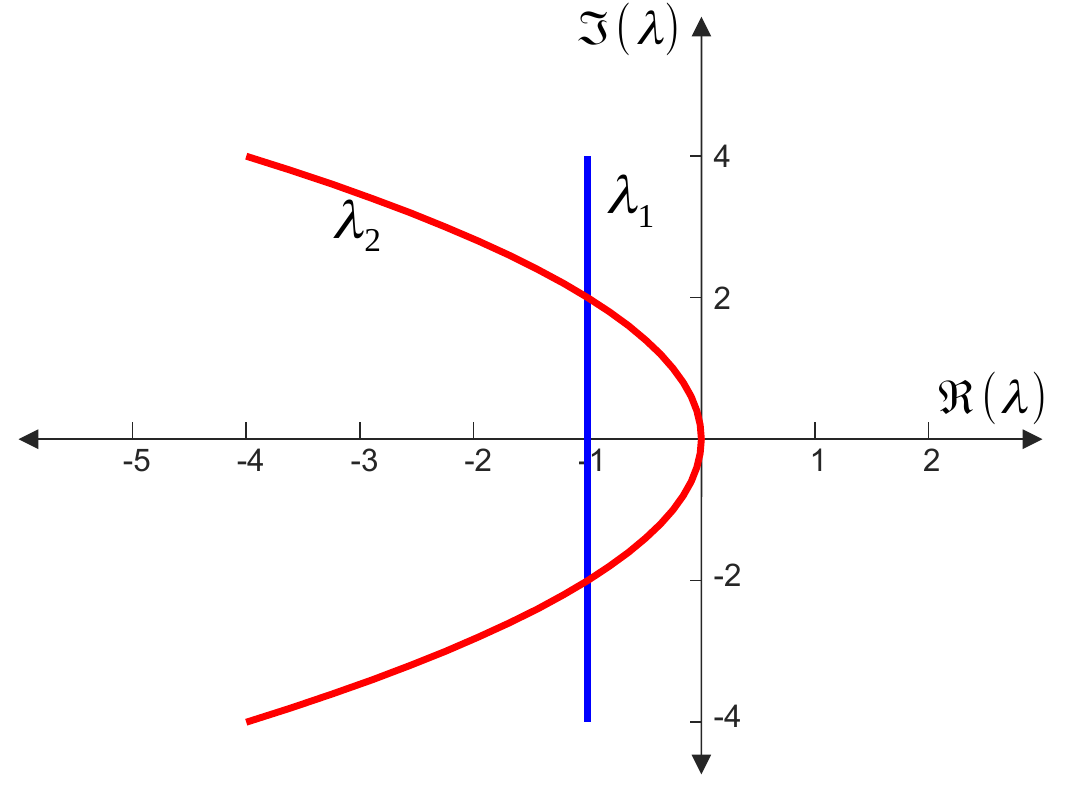}
	\caption{The essential spectrum of $\mathcal{L}$ (for travelling wave S1) is given by $\lambda$ along the Fredholm borders, in this case the union of a straight line ($\lambda_1$) and a parabola ($\lambda_2$) in the complex plane. 
	}
	\label{fig_essential}
\end{figure} 
\subsection{Essential spectrum} 
\RM{The essential spectrum is the spectrum up to relatively compact perturbations.}
To this end we consider the asymptotic matrices $A^\pm(\lambda)=\lim_{z \to \pm \infty} A(z,\lambda)$ given by
\begin{equation}\label{eq6}
A^+(\lambda)=\begin{pmatrix} 0&&\frac{-s^3}{\kappa\,(s-1)^3}&&0\\
    (\frac{1}{s}-1)\,\lambda&& \frac{s^2}{\kappa\,(1-s)}&& 0\\
    0&& \frac{1}{1-s}&& \frac{\lambda+1}{s-1}
    \end{pmatrix}\quad \text{and} \quad
            A^-(\lambda)=\begin{pmatrix} 0&&\frac{-1}{\kappa}&&0\\
    -\lambda&& -\frac{s}{\kappa}&& 0\\
    0&& -\frac{1}{s}&& \frac{\lambda+1}{s}
    \end{pmatrix}.
    \end{equation}   
The asymptotic operator of $\mathcal{T}$ is given by
\begin{equation}
  \mathcal{T}_\infty \begin{pmatrix}
     \delta\rho\\
    \delta v\\
      \delta a
     \end{pmatrix} =\begin{pmatrix}
     \delta\rho'\\
    \delta v'\\
      \delta a'
     \end{pmatrix}-A_\infty\begin{pmatrix}
    \delta\rho\\
    \delta v\\
      \delta a
     \end{pmatrix} \; , \quad 
     \text{where} \quad 
     A_\infty=
         \begin{cases}
         A_-(\lambda)   &z<0 \; ,\\
         A_+(\lambda)   &z\geq 0 \; .
         \end{cases}
\end{equation}
Note that the opterator $\mathcal{L}$ is a relatively compact perturbation of the asymptotic operator $\mathcal{L}_\infty$ which is defined as the limit of $\mathcal{L}$ as $z\to\pm\infty$ and which is equivalent to the operator $\mathcal{T_\infty}$ (\cite{a24}, Theorem 3.1.11), see also \cite{a26}). According to Weyl's Essential Spectrum Theorem (\cite{a24}, Theorem 2.2.6), the operators $\mathcal{L}$ and $\mathcal{L}_\infty$ have the same spectra, i.e., $\sigma_{ess}(\mathcal{L})=\sigma_{ess}(\mathcal{L}_\infty)$, or equivalently $\sigma_{ess}(\mathcal{T})=\sigma_{ess}(\mathcal{T}_\infty)$. 

Exponential dichotomies can be used to characterise the spectrum of an operator. According to this concept each solution to \eqref{equ_linearised} decays exponentially either for $z \rightarrow \infty$ or for $z \rightarrow -\infty$ \cite{a26}.
For spatially constant matrices the presence of an exponential dichotomy implies that the matrix is hyperbolic.
The Morse indices of the constant matrices $A_{\pm}(\lambda)$ are defined as the dimension of their unstable subspaces written as $i_{\pm}(\lambda)$, respectively. For $\lambda\in\mathbb{C}$ such that $\mathcal{T_\infty}$ is Fredholm, we have ind$(\mathcal{T_\infty}-\lambda)=i_-(\lambda)-i_+(\lambda)$ (\cite{a24}, Lemma 3.1.10). As a consequence we can define the essential spectrum of $\mathcal{L_\infty}$ as (\cite{a24})
\begin{equation}
	\sigma_{ess}(\mathcal{L_\infty})=\{\lambda\in\mathbb{C}\mid i_-(\lambda)\neq i_+(\lambda)\}\cup\{\lambda\in\mathbb{C}\mid  dim\; \mathbb{E}^c(A_\pm(\lambda))\neq 0\} \;,
\end{equation}
where $\mathbb{E}^c$ denotes the centre subspace associated to the asymptotic linearised system.

The dispersion relations of $A^-$ and $A^+$ which characterise $ dim\; \mathbb{E}^c(A_\pm(\lambda))\neq 0$ are defined by
\begin{equation}
\det(A^- - i\, \mu  \, \rm{Id})=0 \quad \text{and} \quad \det(A^+-i\, \mu \, \rm{Id})=0\;,\quad \text{where} \quad \mu\in\mathbb{R} \; .
\end{equation}
These relations of both $A^-$ and $A^+$, after rescaling the spatial eigenvalue $\mu$, coincide and are given by
\begin{equation}
\label{disp_rel}
\lambda^{1} =-1+i\mu s \quad \text{and}\quad \lambda^{2}=-\mu^2 \kappa+i\mu s \; .
\end{equation}
Note that the matrix eigenvalues of $A^{-}(\lambda)$ defined in \eqref{eq6} are given by
\begin{equation}
\label{equ_muminus}
\mu_1(\lambda)= \frac{\lambda+1}{s}\; , \quad 
\mu_2(\lambda)=\frac{-s}{2\kappa}+\frac{\sqrt{s^2+4\kappa\lambda}}{2\kappa}\; , \quad
\mu_3(\lambda)=\frac{-s}{2\kappa}- \frac{\sqrt{s^2+4\kappa\lambda}}{2\kappa}\;.
\end{equation}
and those of $A^{+}(\lambda)$ are
\begin{align}
\label{equ_muplus}
\mu_1(\lambda)&= \frac{\lambda+1}{s-1}\quad \text{and}\quad
\mu_{2,3}(\lambda)=\frac{s^2}{2\kappa(1-s)}\mp s\,\frac{\sqrt{s^2+4\kappa\lambda}}{2\kappa(1-s)} \; ,
\end{align}
which are multiples of the matrix eigenvalues of $A^{-}(\lambda)$ by $\frac{s}{s-1}>0$. As a consequence the Morse indices coincide for all $\lambda$, i.e. $i_-(\lambda)=i_+(\lambda)$.
Therefore the essential spectrum only consists of the the Fredholm borders \eqref{disp_rel}
which are visualised in Fig.~\ref{fig_essential}.

\subsection{Absolute spectrum} 
While the absolute spectrum is not spectrum \cite{a27}, it gives information about the stability of the operator $\mathcal{L}$ in  exponentially weighted spaces. Most notably it tells us how far the essential spectrum may be shifted to the left by considering weighted function spaces with exponential weights at $\pm \infty$. If the absolute spectrum lies completely in the open left half-plane of the complex plane, then we say $\mathcal{L}$ is absolutely stable, otherwise absolutely unstable.

We note that common value of the of the Morse index for $\lambda \gg 1$ by $i_{\infty}$, i.e. $i_{\infty}=i_{\pm}(\lambda)$. We define the matrix (spatial) eigenvalues $\mu_{\pm}^{j}\,,\,j=1,...,n\,$ of the asymptotic matrices $A_{\pm}(\lambda)$ ordered according to the size of their real parts\cite{a24},
\begin{equation*}
    Re\,\mu_{\pm}^{1}(\lambda)\geq.....\geq\,Re\,\mu_{\pm}^{l}(\lambda)\geq.....\geq\,Re\,\mu_{\pm}^{n}(\lambda)\; ,
\end{equation*}
and introduce the stable and unstable extrema 
\begin{equation*}
    \mu_{\pm}^{u}(\lambda)=Re\,\mu_{\pm}^{i_\infty}(\lambda) \; , \quad  \mu_{\pm}^{s}(\lambda)=Re\,\mu_{\pm}^{i_\infty+1}(\lambda) \; .
\end{equation*}
Hence $\mu_{\pm}^{u}(\lambda) $ denotes the smallest (positive) real part of any of the matrix eigenvalues, and  $\mu_{\pm}^{s}(\lambda) $ denotes the largest(negative) real part of any of the matrix eigenvalues, of the asymptotic matrices $A_{\pm}(\lambda)$. Here we have $\mu_{\pm}^{u}(\lambda)>0>\mu_{\pm}^{s}(\lambda)$ for $Re\,\lambda\gg1$. As we move the eigenvalue $\lambda$ towards the Fredholm border coming from the very far right of the complex plane, then at least one of the matrix eigenvalues will become close to the imaginary axis. Moreover, for the weighted spaces, when the distance between the spatial eigenvalues, $\mu_{\pm}^{u}\,(\lambda)-\mu_{\pm}^{s}\,(\lambda)$, becomes zero, we can not choose a weight that renders the operator $\mathcal{L}-\lambda$ Fredholm with index zero. This motivates the following definition.

The subset $\sum_{abs}^{+}$\, of $\mathbb{C}$ consists exactly of those $\lambda$ for which $Re\,\mu_{+}^{i_\infty}(\lambda)\,=\,Re\,\mu_{+}^{i_{\infty}+1}(\lambda)\,$. Analogously, $\lambda$ is in $\sum_{abs}^{-}$ if, and only if, $Re\,\mu_{-}^{i_\infty}(\lambda)\,=\,Re\,\mu_{-}^{i_{\infty}+1}(\lambda)\,$. Finally, we say that $\lambda$ is in the absolute spectrum $\sum_{abs}$ of an operator if $\lambda$ is in $\sum_{abs}^{+}$ or in $\sum_{abs}^{-}$ (or in both).

\begin{figure}
	\centering
	\includegraphics[width=0.8\textwidth]{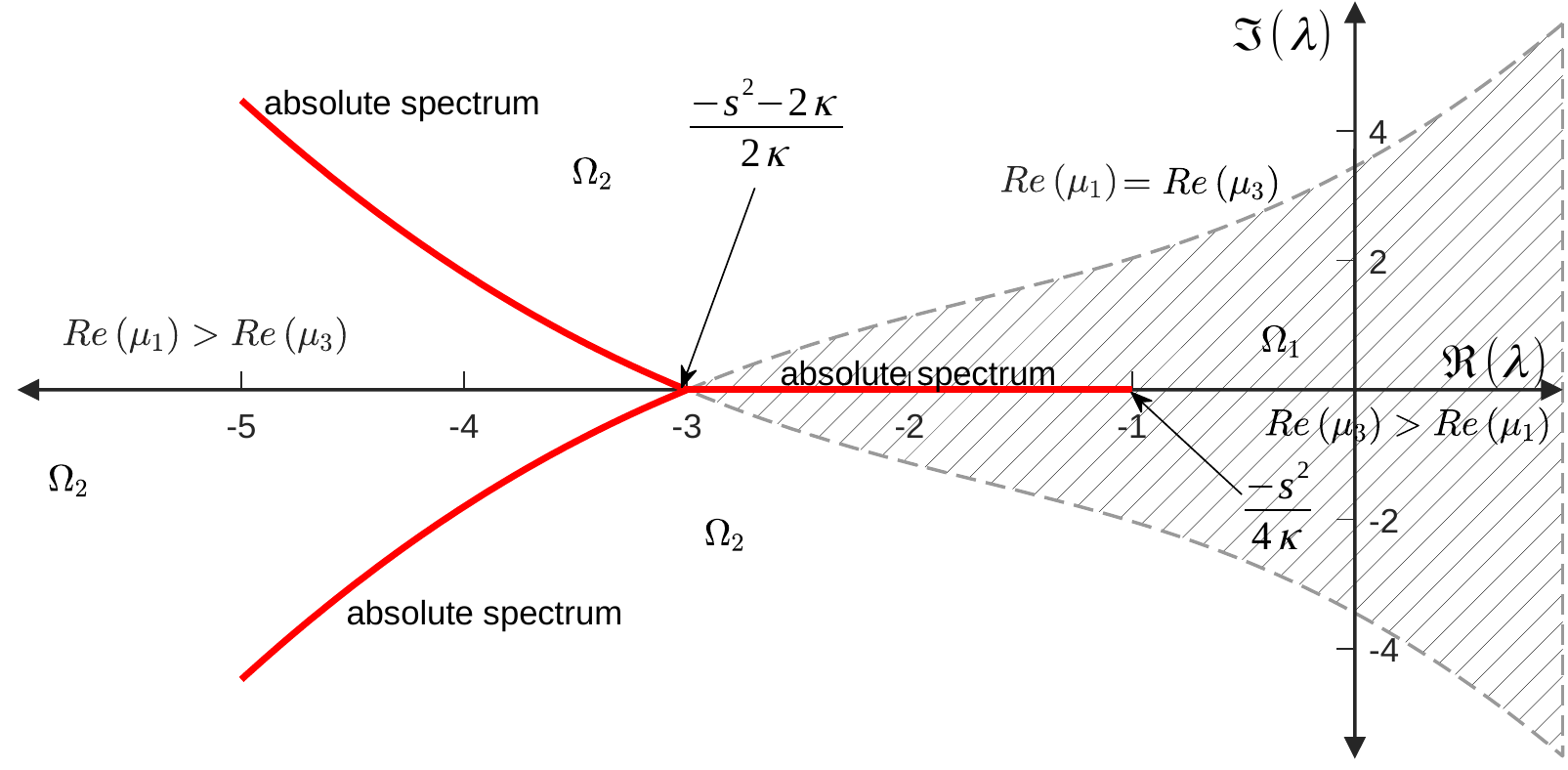}
	\caption{Absolute spectrum shown in red (for $s=-2$, $\kappa=1$).}
	\label{fig_absspectrum}
\end{figure}

For $ Re(\lambda) \gg 1$ the matrix eigenvalues of $A^{-}(\lambda)$ given in \eqref{equ_muminus} satisfy $Re(\mu_2)>0$ as well as $Re(\mu_3)<0$ and $Re(\mu_1)<0$ (since $s<0$). 
In order to find which of the two asymptotically negative spatial eigenvalues is larger (real part) for a given $\lambda$, we solve the inequality $Re(\mu_1)\leq Re(\mu_3)$, i.e., 
\begin{align*}
  Re\left(\frac{1+\lambda_1+i\lambda_2}{s}\right) &\leq Re \left( \frac{-s}{2\kappa}-\frac{\sqrt{s^2+4\kappa(\lambda_1+i\lambda_2)}}{2\kappa}\right).
 \end{align*}

 In this case we get, writing $S(\lambda_1)=\frac{ (s^2+2\kappa\,(1+\lambda_1))\,\sqrt{(s^2+\kappa\,(1+\lambda_1)^2)}}{s^2\sqrt{\kappa}}$, that
 \begin{equation*}
 \Omega_1=\left\{\lambda=\lambda_1+i\,\lambda_2 \,\,\text {such that}\,\, \lambda_1\in \left[-\frac{\,s^2+2\kappa}{2 \kappa},\infty\right)\,\, \text{and} \,\,  \lambda_2\in \left [- S(\lambda_1), S(\lambda_1) \right]\right\}
 \end{equation*}
 (hatched area in Fig.~\ref{fig_absspectrum}) and $Re(\mu_1) > Re(\mu_3)$ is satisfied by all $\lambda$ in the complement $\Omega_2= (\mathbb{C} \setminus \Omega_1)$ of this set.
As a consequence the absolute spectrum is the set of all $\lambda$ such that $Re(\mu_2)=Re(\mu_3)$ in $\Omega_1$ and $Re(\mu_2)=Re(\mu_1)$ in $\Omega_2$.

Now for the absolute spectrum in the set $\Omega_1$, we solve the equation $Re(\mu_2)=Re(\mu_3)$, i.e.
\begin{align*}
  Re \left( \frac{-s}{2\kappa}+\frac{\sqrt{s^2+4\kappa\lambda}}{2\kappa}\right)&=Re\left(\frac{-s}{2\kappa}-\frac{\sqrt{s^2+4\kappa\lambda}}{2\kappa}\right),
\end{align*}
which implies that $Re(\lambda)\in\left[-\frac{s^2+2\kappa}{2\kappa},-\frac{s^2}{4\kappa}\right]$ and $Im(\lambda)=0$.
\par
To obtain the absolute spectrum in the set $\Omega_2$, we solve the equation $Re(\mu_1)=Re(\mu_3)$, i.e.
\begin{align*}
 Re\left(\frac{1+\lambda_1+i\lambda_2}{s}\right) &=Re \left( \frac{-s}{2\kappa}+\frac{\sqrt{s^2+4\kappa(\lambda_1+i\lambda_2)}}{2\kappa}\right), 
\end{align*}
which is the case for all $\lambda=\lambda_1+i\,\lambda_2$ such that $\lambda_1\,<-\frac{\,s^2+2\kappa}{2\kappa}\,\, \text{and}\,\, \lambda_2=  \pm \frac{ (s^2+2\kappa\,(1+\lambda_1))\,\sqrt{(s^2+\kappa\,(1+\lambda_1)^2)}}{s^2\sqrt{\kappa}} $. 
Thus the absolute spectrum of $A^{-}(\lambda)$ is given by
\begin{multline}
\label{eq11}
\sigma_{abs}^{-}=\left[-\frac{s^2+2\kappa}{2\kappa},-\frac{s^2}{4\kappa}\right] \,\cup \Bigg\{\lambda=\lambda_1+i\,\lambda_2 \,\mid  \lambda_1\,<-\frac{\,s^2+2\kappa}{2\kappa}\,\, \text{and} \\
 \lambda_2=  \pm \frac{ (s^2+2\kappa\,(1+\lambda_1))\,\sqrt{(s^2+\kappa\,(1+\lambda_1)^2)}}{s^2\sqrt{\kappa}} \Bigg\}.
\end{multline}

Finally the matrix eigenvalues of $A^{+}(\lambda)$ given in \eqref{equ_muplus} are a multiple of the matrix eigenvalues of $A^{-}(\lambda)$, namely by $\frac{s}{s-1}>0$. So the absolute spectrum of $A^{+}(\lambda)$ is equal to the absolute spectrum of $A^{-}(\lambda)$ and the absolute spectrum is given by \eqref{eq11} (see Fig.~\ref{fig_absspectrum}).
 
 \RM{This implies in a conveniently chosen weighted space, and without the presence of any point spectrum with non-negative real part, save for at $\lambda =0$ (see Section \ref{subsec:pointspec}), the wave is spectrally stable. To determine the weights for which we potentially have spectral stability, we follow \cite{DvHM17,a27} looking for a so-called {\em ideal weight}. The ideal weight will be the (two-sided) weight which maximises the resolvent set.
We find the ideal weight $\eta^{*}_{+}$ to satisfy $\textrm{Re}(\eta^*_{+}) = \frac{s^2}{2\kappa(1-s)}$, for the operator corresponding to $A^+(\lambda)$, and $\textrm{Re}(\eta^*_{-})=\frac{-s}{2 \kappa}$ for the operator corresponding to $A^-(\lambda)$. We note that these also differ by a factor of $\frac{s}{s-1}$.  By considering perturbations that decay at least like $\exp{(-\eta z)}$ as $z \to \infty$ with $\eta \in (0,\frac{s^2}{\kappa(1-s)}]$, and like $\exp{\eta z}$ as $z \to -\infty$, with $\eta \in (0,-\frac{s}{\kappa}]$, we have that the essential spectrum will be contained in the (open) left half plane. This is worth noting because the derivative of the S1 wave decays as $z \to +\infty$ like $\exp{\left(\frac{z}{s-1}\right)}$, and as $z\to -\infty$, like $\exp{\left(\frac{-3sz}{\kappa}\right)}$ and so will remain in the weighted space for the weights we want to consider. This means that $\lambda = 0$ will still be an eigenvalue both in $L^2$ and in the weighted space. As we shall see in the next section, $\lambda = 0$ is the only element of the point spectrum that we can numerically find. 
 }

\subsection{Point Spectrum}\label{subsec:pointspec}
In this section we investigate the point spectrum of $\mathcal{L}$. We compute the Evans function whose zeros in the complex plan characterise the point spectrum. 
In 1972, J. W. Evans used this technique to investigate the stability of the solution to equations modelling the nerve axon \cite{a22,a23}. He showed that the Evans function $D(\lambda)$ is always analytic to the right side of the essential spectrum. In addition it has an analytic extension up to the absolute spectrum (Fig.~\ref{fig_absspectrum}).

To compute the Evans function, we start by identifying the stable and unstable eigenvectors of the asymptotic matrices $A^{\pm}(\lambda)$ respectively in an attempt to construct an integrable function which plays the role of an eigenvector for the linearised system. 

The unstable eigenvalue of $A^{-}$ and the stable eigenvalues of $A^{+}$  are given by
\begin{equation}\label{eq12}
\mu_2^{-}(\lambda)=\frac{-s+\sqrt{s^2+4\kappa \lambda}}{2\kappa} \; , \quad 
\mu_1^{+}(\lambda)= \frac{\lambda+1}{s-1} \quad \text{and} \quad 
\mu_3^{+}(\lambda)= \frac{s^2+s\sqrt{s^2+4\kappa \lambda}}{2\kappa(1-s)}  \;,
\end{equation}
which satisfy $\mu_2^{-}(\lambda)>0$ and $\mu_1^{+}(\lambda)\,,\mu_3^{+}(\lambda)<0$ for ${\rm Re} \, \lambda \gg 1$ since that travelling wave speed is negative, $s=s_1<0$. 

We introduce the eigenvectors $X_0^-(\lambda)$ associated to the eigenvalue $\mu_2^{-}(\lambda)$ of $A^-$ as well as $X_0^+(\lambda)$ and $Y_0^+(\lambda)$ being eigenvectors of $A^+$ associated to $\mu_1^{+}(\lambda)$ and $\mu_3^{+}(\lambda)$.

To construct an eigenfunction of the linearised system we solve \eqref{final_A} using these vectors as initial, respectively terminal conditions until $z=0$. Then
the Evans function is defined as the Wronskian
\begin{equation}\label{equ_evans}
D(\lambda)={\rm det}\; [X^{-}(z=0,\lambda),X^{+}(z=0,\lambda), Y^{+}(z=0,\lambda)] \; ,
\end{equation}
which vanishes if the stable and unstable eigenvectors propagated to $z=0$ are linearly dependent and can be combined into a smooth eigenfunction.

We start by computing $X^+(z=0, \lambda)$ and $Y^+(z=0, \lambda)$ which can be done in terms of a closed-form expression. Only for the computation of $X^-(z=0, \lambda)$ we will resort to numerical results.

\bigskip 

For $z>0$ it holds that $R  \equiv s/(s-1)$ and the matrix $A(z,\lambda)$ from \eqref{final_A} is given by
\begin{equation}
\label{eqevan1}
A(z,\lambda)=\begin{pmatrix} 0&&\frac{-s^3}{\kappa \,(s-1)^3}&&0\\
  (\frac{1}{s}-1)\,\lambda&& \frac{s^2}{\kappa \,(1-s)}&& 0\\
  0&& \frac{A'-1}{s-1}&& \frac{\lambda+1}{s-1}
  \end{pmatrix}\,.
\end{equation}
Note that due to the zeros in the third column the equations for $\delta \rho$ and $\delta v$ are not coupled to $\delta a$. They satisfy
\begin{equation}\label{eqevan2}
\begin{pmatrix} \delta \rho' \\ \delta v' \end{pmatrix}=
\begin{pmatrix} 0&&\frac{-s^3}{\kappa \,(s-1)^3}\\
(\frac{1}{s}-1)\,\lambda&& \frac{s^2}{\kappa \,(1-s)}
\end{pmatrix}
\begin{pmatrix} \delta \rho \\ \delta v \end{pmatrix} \; ,
\end{equation}
which is a system of two linear, constant-coefficient equations. It admits one fundamental solution with negative eigenvalue (we omit the unstable fundamental solution) given by 
\begin{equation}\label{eqevan3}
\begin{pmatrix}
\delta\rho\\
\delta v
\end{pmatrix} = C_{3}\,\exp(\mu_{3}^+ (\lambda) \,z ) \begin{pmatrix}
\frac{s^2\,(-s+\sqrt{s^2+4\,\kappa \,\lambda})}{(s-1)^2\,\kappa}\\
2\,\lambda
\end{pmatrix}
\end{equation}

With this information we can rewrite the equation for $\delta a$ which is contained in
the third row of \eqref{eqevan1},
\begin{align*}
   \delta a '&=\frac{A'-1}{s-1}\,\delta v+\frac{\lambda+1}{s-1}\, \delta a\\
    &=\frac{A'-1}{s-1}\,2\,\lambda\,C_{3}\,e^{\mu_{3}\,z} +\mu_1^+\, \delta a \; .
\end{align*}
The general solution is given by 
\begin{displaymath}
\delta a = C_1 e^{\mu_1^+ z}+ \frac{2 \lambda}{s-1} C_3 \int_0^z e^{(\mu_3^+-\mu_1^+) z}  \left(\frac{\alpha-1}{s-1}e^{\frac{z}{s-1}}-1 \right) \, d\tilde z
\end{displaymath}

This implies that 
$$
\lim_{z\to 0^{+}} \begin{pmatrix}
    \delta\rho\\
    \delta v\\
     \delta a
     \end{pmatrix}=\begin{pmatrix}
C_3\,\frac{s^2\,(-s+\sqrt{s^2+4\,\kappa\,\lambda})}{(s-1)^2\,\kappa}\\C_3\,
2\,\lambda\\
C_1
\end{pmatrix} \; .
$$ 
While most components of $A(z,\lambda)$ in \eqref{final_A} are functions, the derivative of the active speed of migration \eqref{equ_M} is a $\delta$-distribution, $M'(a)=\delta(a-\alpha)$. The change of variables between $z$ and the monotone function $a=A(z)=A_1(z)$ shows that when integrating with respect to $z$ the following expression is a $\delta$-distribution centred at $z=0$, $M'(A(z)) A'(z)=\delta (z)$. Therefore $\delta \rho$ and $\delta v$ in the solution of \eqref{final_A} undergo a jump at $z=0$ which involves the factor $1/A'(0)$ (for details see appendix~\ref{sec_jump}). The left limit as $z\rightarrow 0$ of the solution vector is given by 
 \begin{align*}
 \lim_{z\to 0^{-}} \begin{pmatrix}
     \delta\rho\\
    \delta v\\
      \delta a
     \end{pmatrix}=&\begin{pmatrix}
C_3\,\frac{s^2\,(-s+\sqrt{s^2+4\,\kappa \,\lambda})}{(s-1)^2\, \kappa }\\C_3\,
2\,\lambda \\
 C_1
 \end{pmatrix}-
\frac{C_1}{A'(0)} \frac{1}{\kappa} \begin{pmatrix}
R^3(0)\\ R(0)\,s\\
 0
 \end{pmatrix}
 \\
 =& C_3\,\begin{pmatrix}
\frac{s^2\,(-s+\sqrt{s^2+4\,\kappa \,\lambda})}{(s-1)^2\,\kappa}\\
2\,\lambda\\
0
 \end{pmatrix}+\,C_1\,\begin{pmatrix}
 -\frac{s^3}{(s-1)^2\,(\alpha-1)\,\kappa}\\
-\frac{s^2}{(\alpha-1)\,\kappa}\\ 1
 \end{pmatrix}.
  \end{align*}
where we take the values of $R(0)$ and $A'(0)$ from the travelling wave profile S1.

We obtain the closed-form solutions $X^{+}(0,\lambda)$ and $Y^{+}(0,\lambda)$ corresponding to the stable eigenvalues $\mu_1^{+}(\lambda)$ and $\mu_3^{+}(\lambda)$, 
\begin{align}
\label{equ_XY}
    X^{+}(0,\lambda) =\begin{pmatrix}
 -\frac{s^3}{(s-1)^2\,(\alpha-1)\,k}\\
-\frac{s^2}{(\alpha-1)\,k}\\ 1
 \end{pmatrix} 
 \quad \text{and}\quad
 Y^{+}(0,\lambda)=\begin{pmatrix}
\frac{s^2\,(-s+\sqrt{s^2+4\,k\,\lambda})}{(s-1)^2\,k}\\
2\,\lambda\\
0
 \end{pmatrix}.
\end{align}

Finally, in order to evaluate the Evans function \eqref{equ_evans} for a given $\lambda \in \mathbb{C}$, we compute $X^-(0,\lambda)$ propagating the eigenvector associated to the spatial eigenvalue $\mu_2$ given by 
$$
v_2=\begin{pmatrix}
\frac{-s^2+\sqrt{s^4+4s^2\kappa\lambda}}{2\lambda s \kappa}\\
1\\ \frac{2\kappa}{s^2+2\kappa(1+\lambda)+\sqrt{s^4+4s^2\kappa\lambda}}
\end{pmatrix} \; \text{for} \; \lambda \neq 0 
\quad \text{and} \quad 
v_2=
\begin{pmatrix}
\frac{1}{s}\\
1\\ \frac{\kappa}{s^2+\kappa} 
\end{pmatrix}  \; \text{for} \; \lambda=0
$$
from an arbitrary small value of $z$ (we choose $z=-20$) until $z=0$ solving the system \eqref{final_A} numerically.

\bigskip 

We apply the Argument principle to identify zeros of the Evans function \cite{a26}. The argument principle states that
\begin{equation*}
    W=\frac{1}{2 \pi i}\,\oint_C\,\frac{D'(\lambda)}{D(\lambda)}\,d\lambda 
\end{equation*}
corresponds to the winding number around the origin of the image of $C$ under the map $D$. Here $C$ is a closed curve, oriented in a counterclockwise direction. Since we choose the contour $C$ to the right of the absolution spectrum the Evans function $D(\lambda)$ is analytic on and inside $C$ and the the winding number $W$ corresponds to the number of zeros of $D(\lambda)$ inside $C$.  
Note that $\lambda=0$ is always an eigenvalue corresponding to the propagation of the travelling wave profile \cite{SANDSTEDE2002}. 

\begin{figure}[H]
   \centering
   \includegraphics[width=.50\textheight  ]{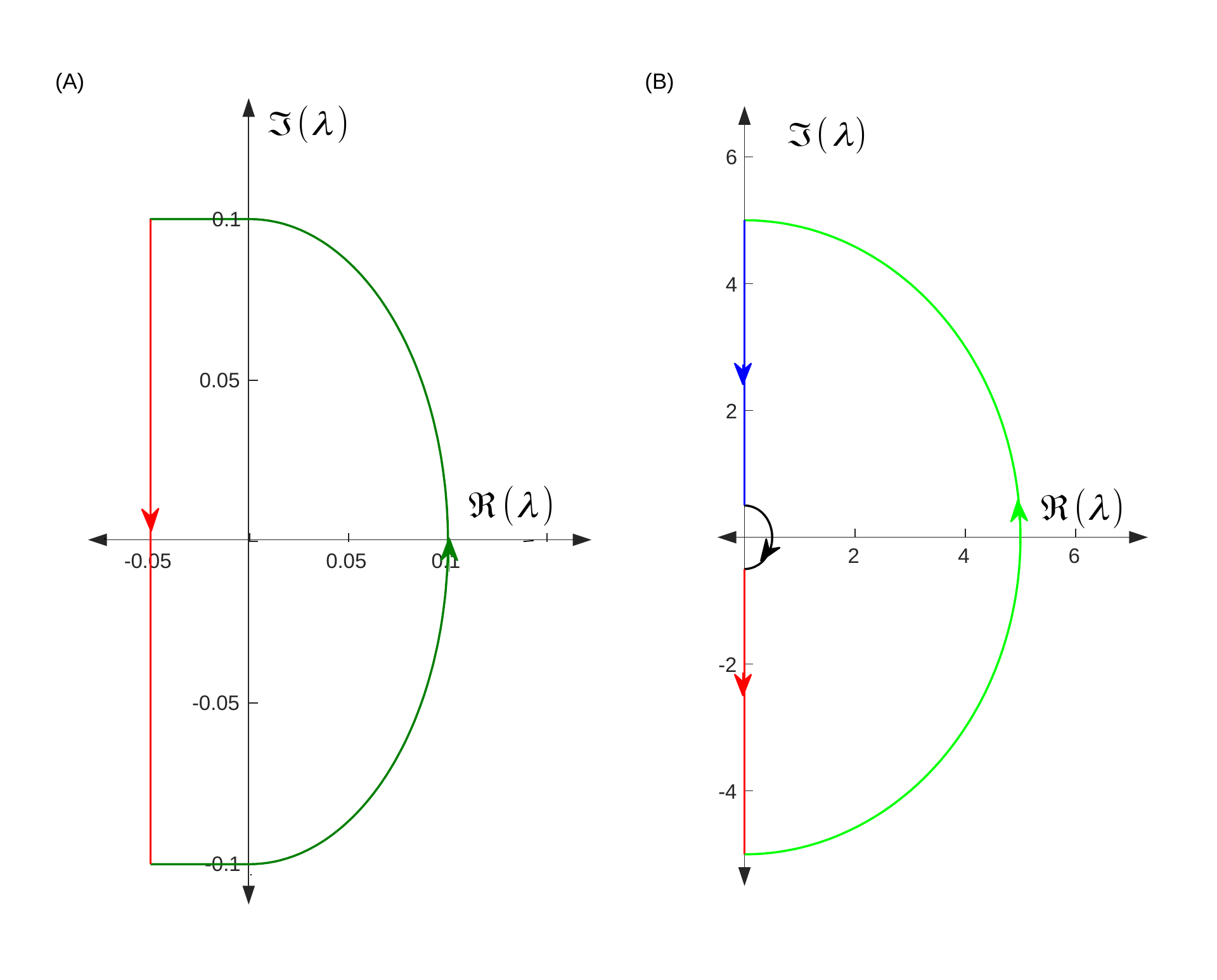}
   \caption{(A) and (B) show the contours $C_1$ and $C_2$ including and excluding the origin, respectively.}
   \label{contour}
\end{figure}
    
We take the contours $C_1$ and $C_2$ such that these are always \RM{to the right of the absolute spectrum}. For the contour $C_1$, we take the line $l$, parallel to the imaginary axis and right to the branch point with distance $d_l$ from the axis and the semicircle right to the origin has radius $r$. In the right hand side of the complex plane, we draw the contour $C_2$ with the radii $r_i$ and $r_o$ of the inner and outer semicircles respectively excluding the origin.

\begin{figure}
\centering
\includegraphics[width=.55\textheight]{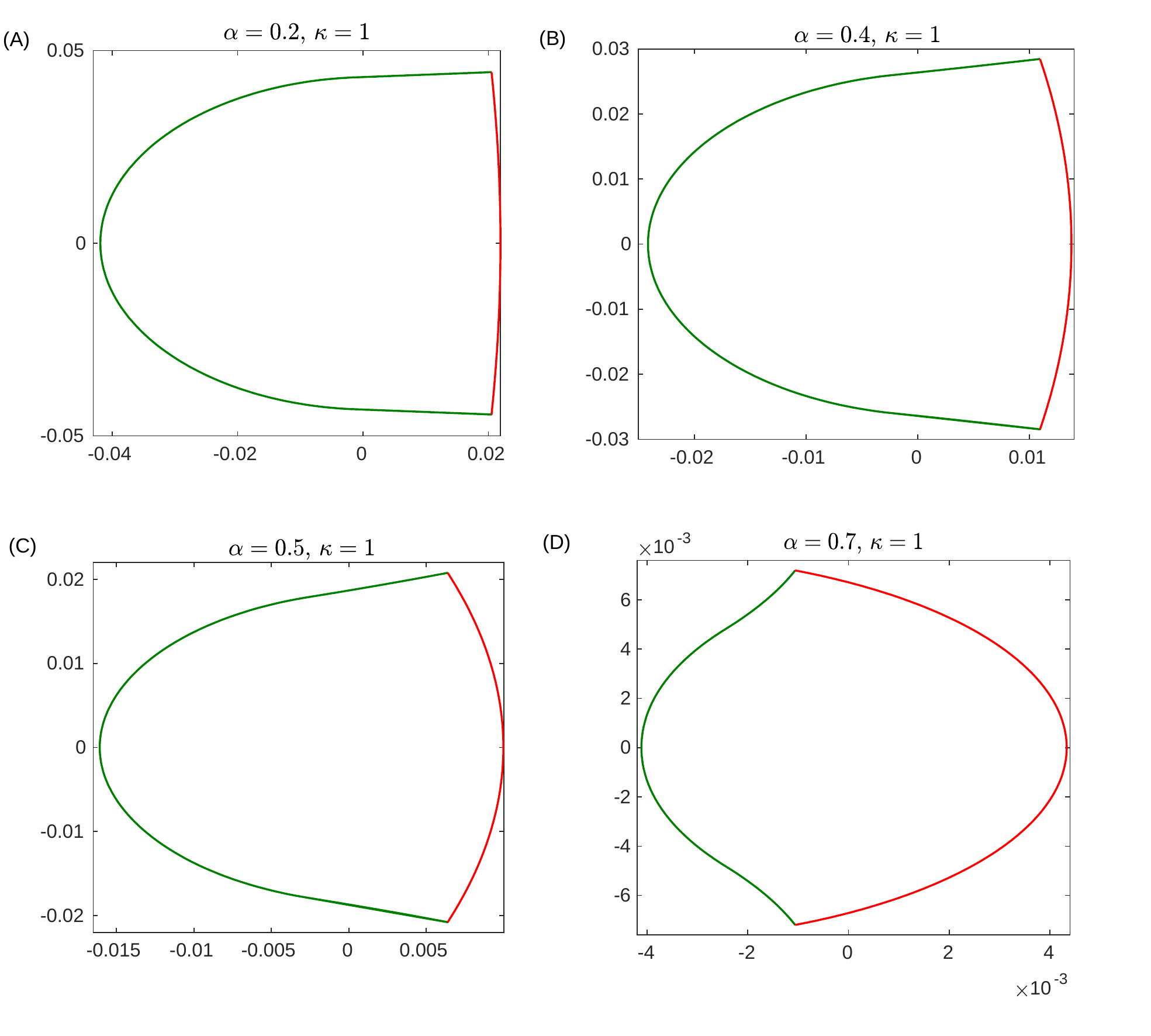}
\caption{$D(C_1)$ for $\alpha=0.2\,,0.4\,,0.5 \,\text{and}\, 0.7$.}
\label{zeros}
\end{figure}
    
In figure \ref{zeros}, we draw the image of $C_1$ (where $d_l=-0.05$ and $r=0.1$) under $D(\lambda)$. According to the argument principle, the winding number for different $\alpha$'s is one, so there is only one zero which is at $\lambda=0$. Hence we see the $0$ is the only eigenvalue close to the origin. 

Now we investigate whether there are eigenvalues on the right side of the complex plane. To do this we define a semicircle excluding the origin. We draw the contour $C_2$ (see Fig.~\ref{contour} (B)). 
\begin{figure}
\centering
\includegraphics[width=.75\textwidth]{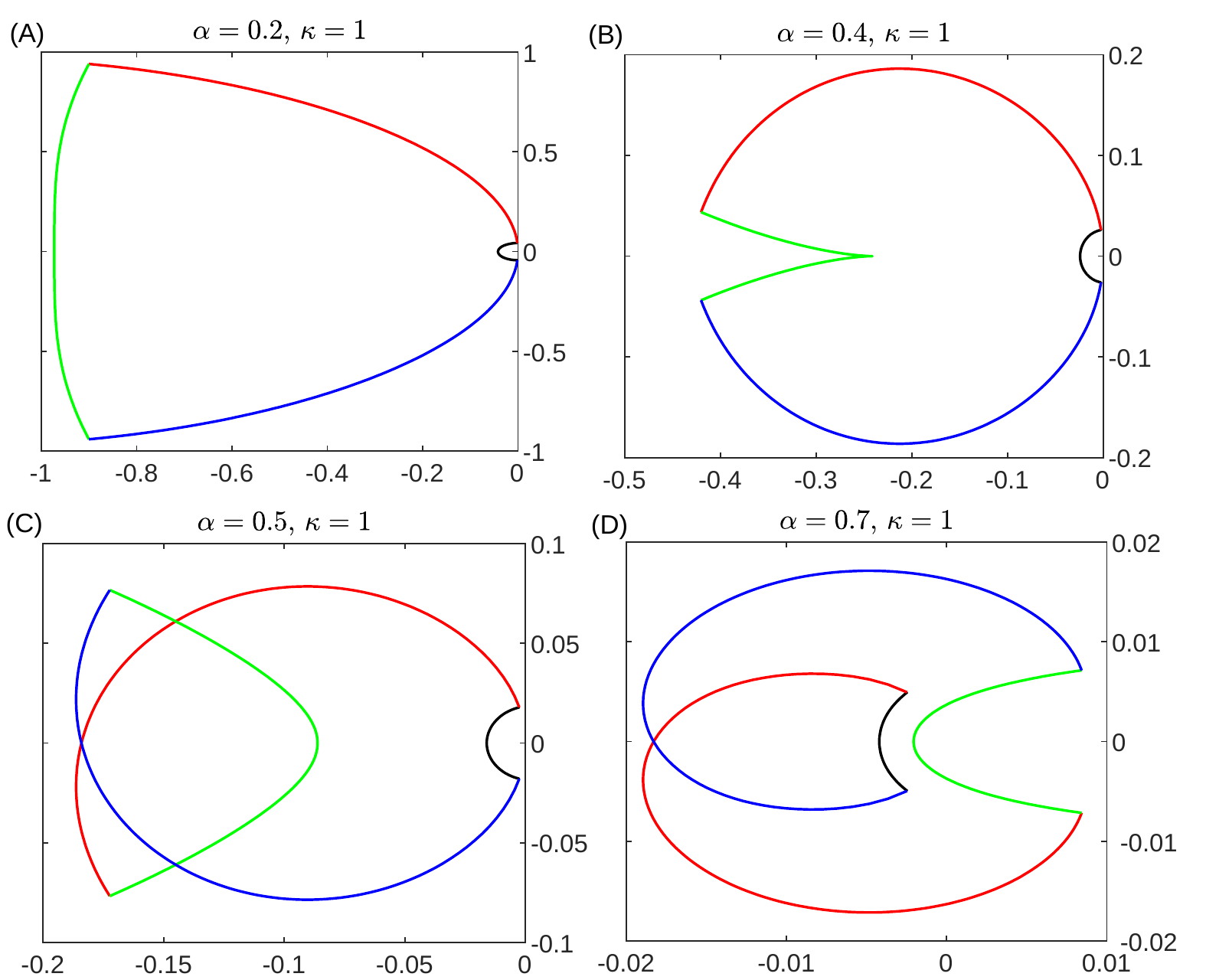}
\caption{$D(C_2)$ for $\alpha=0.2\,, 0.4\,,0.5$ and $0.7$.}
\label{nonzero}
\end{figure}

In figure \ref{nonzero}, we see that the winding number about the origin is zero for the contour $C_2$ with $r_i=0.1$ and $r_o=5$. 

\begin{figure}[H]
	\centering
	\includegraphics[width=.75\textwidth]{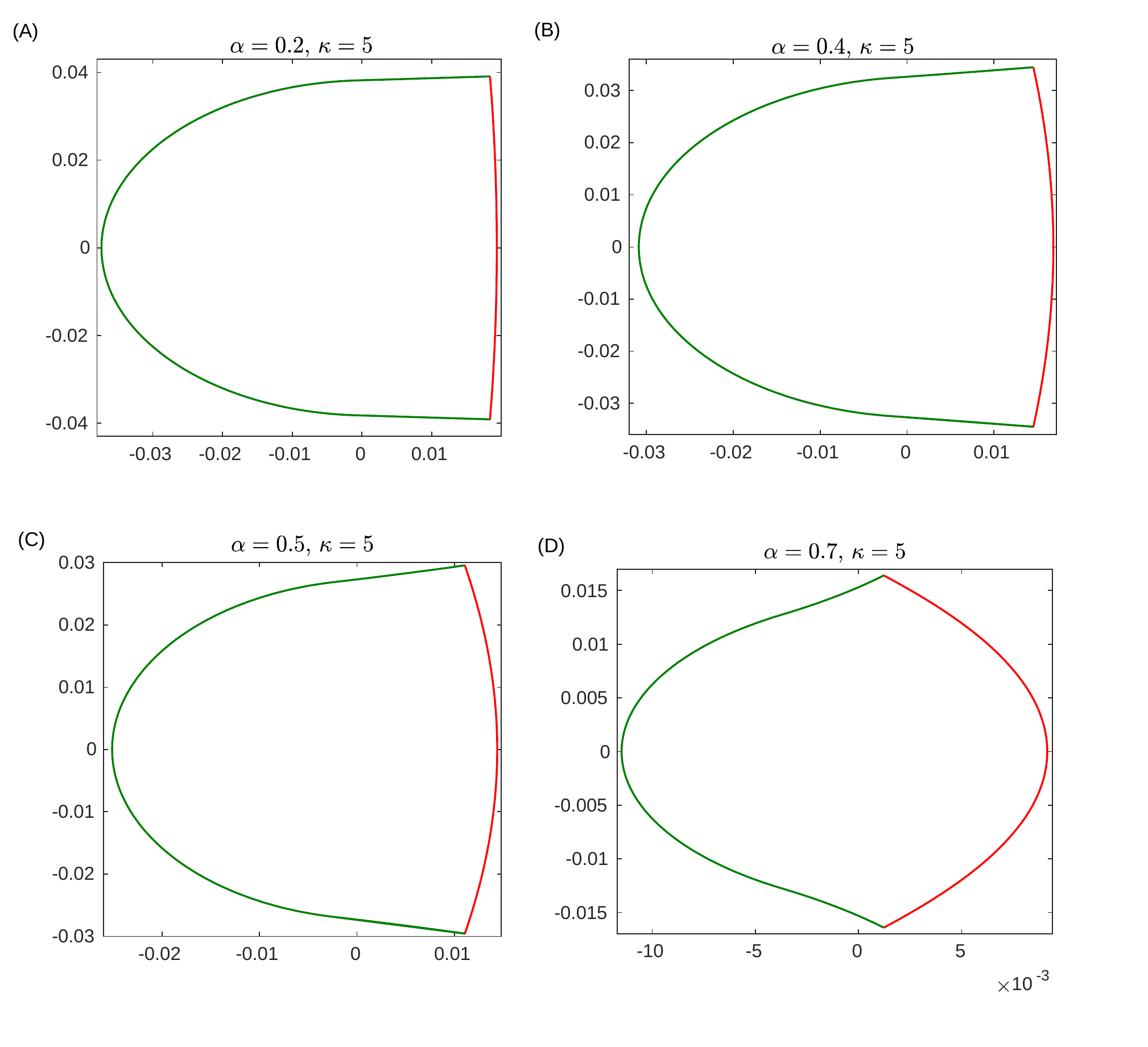}
	\caption{$D(C_1)$ for $\kappa=5$.}
	\label{zero_2}
\end{figure}

\begin{figure}[H]
\centering
\includegraphics[width=.75\textwidth]{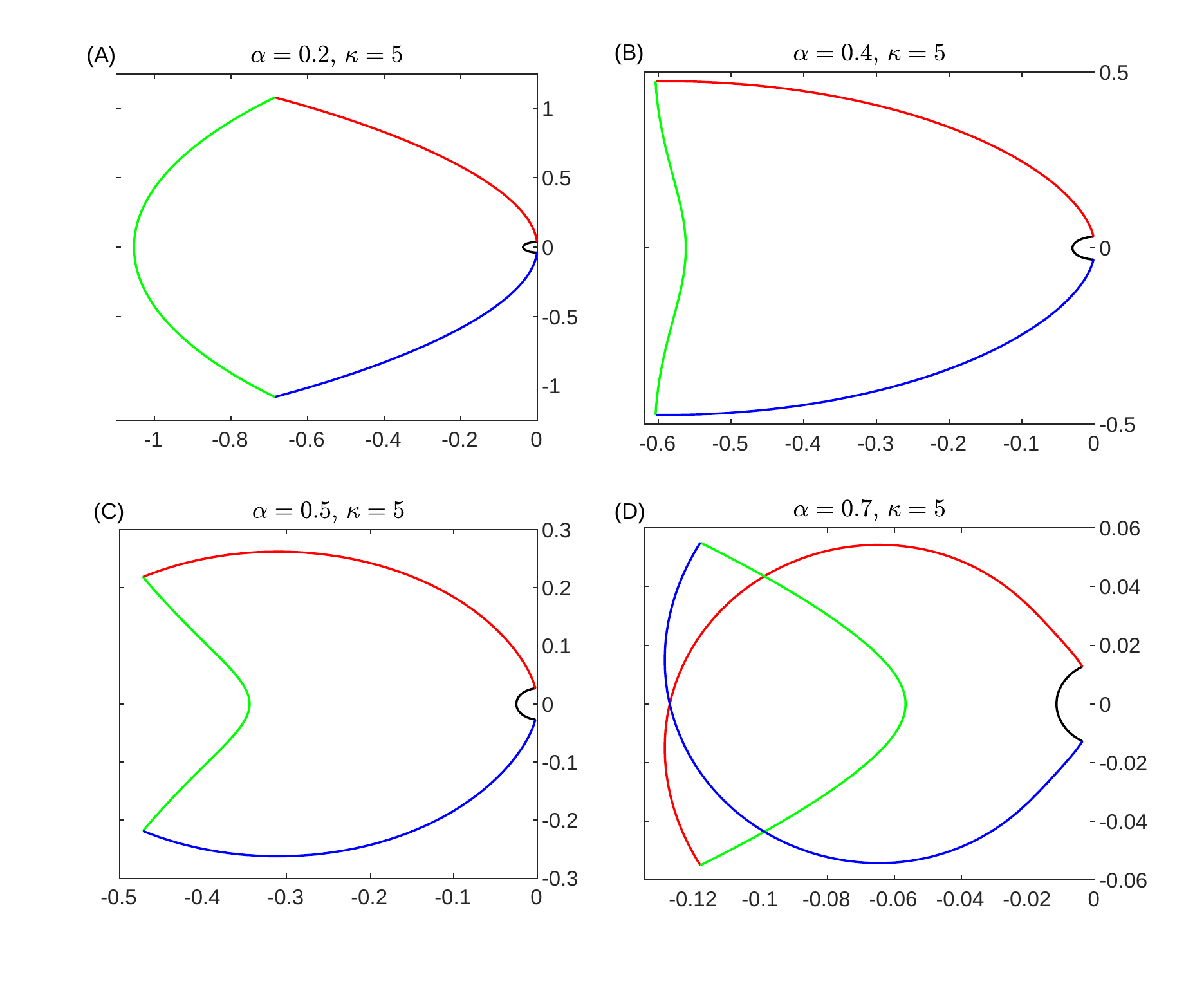}
\caption{$D(C_2)$ for $\kappa=5$.}
\label{nonzero_2}
\end{figure}

This illustrates that there is no point spectrum on the right side of the complex plane other than the eigenvalue at $\lambda=0$ which corresponds to the translocation of the travelling wave. This indicates that the travelling wave solution S1 is spectrally stable. 

\section{Conclusion}
\label{sec_concl}
We study the travelling wave solutions of a 1D model for collective cell migration in an epithelial layer. We identify four different travelling waves, two of them corresponding to gradual polarisation of the cells and the other two corresponding to depolarisation. These travelling wave solutions are related by the two different transformations $T1$ and $T2$.

We apply the transformation $T1$ obtaining the polarisation, respectively depolarisation waves associated to departing cell sheets from the (de)polarisation waves of colliding cell sheets and vice versa. 
The transformation $T2$ can be used to turn polarisation waves into depolarisation waves and vice versa. 

A preliminary test of the stability of these travelling waves computing numerical solutions of the underlaying PDE suggest that they are stable for all parameter values. Yet, we also illustrate that the threshold polarity $\alpha$ which represents the sensitivity to polarisation has a significant impact on the respective sizes of the domains of convergence of the polarisation, respectively depolarisation waves.

Using spectral theory and the Evans function we investigate the stability of travelling waves. This analysis confirms the spectral stability of the travelling wave solutions. While the essential spectrum touches the origin, we show that the absolution spectrum is on the left hand side of the complex plane. As a consequence in appropriately chosen weighted spaces the essential spectrum is contained in the left hand side of the complex plane. The evaluation of the Evans function is done in a way which combines explicit computations and numerics, and illustrates that the only eigenvalue is the simple eigenvalue $\lambda=0$ associated with the propagation of the travelling wave.

We perform this analysis in the context of the polarisation wave associated to departing cell sheets. Due to the smoothness of the map $T_2$ the stability of the polarisation wave also translates into stability of the associated depolarisation wave S3. We add the analogous analysis for the polarisation wave associated to colliding cell sheets in the supplementary material, which through the map $T_2$ also applies to S4.

In conclusion, we find that the travelling waves of cell polarisation (or depolarisation) arising in cell migration due to departing or colliding cell sheets are always stable. For biological tissues this implies that for both, sensitivity to polarisation ($\alpha$) and strength of intercellular mechanical interaction ($\kappa$), there is no absolute threshold such that gradual recruitment of cells into migration, respectively the gradual transition into a non-migratory rest state comes to a halt. Yet, the preliminary numerical experiments (Fig.~\ref{small_alpha}) indicate that the sensitivity to polarisation $\alpha$ characterises behaviour which is reminiscent of a domain of attraction, i.e. lower sensitivity to polarisation (high $\alpha$) renders polarisation waves more prone to disruption.

\section*{Acknowledgments}
DO was supported by ARC Discovery Project DP180102956. RM was supported by ARC Discovery Project DP200102130. NR was supported by an RTP scholarship funded by the University of Queensland (UQ). 

The authors are very grateful to Zoltan Neufeld (UQ) and Hamid Khataee (UQ) for useful conversations and suggestions. 

\appendix
\section{Solution of \eqref{final_A} around $z=0$.}
\label{sec_jump}
The first component of Eq. \eqref{final_A} is given by
\begin{align*}
\delta \rho^{'}&= \frac{3R^{'}}{R}\delta\rho-\frac{R^3}{k}\delta v+ \frac{M'(A)R^3}{k}\delta a \; ,
\end{align*} 
which we write as the 1st order ODE
\begin{align*}
\delta \rho'-C(z)\delta\rho&=D(z) \; ,
\end{align*}
where $C(z)=\frac{3R^{'}}{R}$ and $D(z)=-\frac{R^3}{k}\delta v+ \frac{M'(A)R^3}{k}\delta a$.
Note that the definition of $M=M(A)$ in \eqref{equ_M} implies that $M'(A)=\delta(A-\alpha)$. 
Multiplying both sides by the integrating factor $e^{-\int_{0}^{z} C(\tilde{z})d\tilde{z}}$ we find
\begin{align*}
(e^{-\int_{0}^{z} C(\tilde{z}) d\tilde{z}}\,\delta\rho)'&=D(z)\, e^{-\int_{0}^{z} C(\tilde{z}) d\tilde{z}} \; .
\end{align*}
Integrating on $[-\varepsilon$, $\varepsilon]$ for a small value $\varepsilon>0$ we obtain
\begin{align*}
e^{-\int_{0}^{\varepsilon} C(\tilde{z})
	d\tilde{z}}\,\delta\rho(\varepsilon)-e^{-\int_{0}^{-\varepsilon} C(\tilde{z})
	d\tilde{z}}\,\delta\rho(-\varepsilon)&=\int_{-\varepsilon}^{\varepsilon} D(\hat{z})\, e^{-\int_{0}^{\hat{z}} C(\tilde{z}) d\tilde{z}} d\hat{z}\\
&=\int_{-\varepsilon}^{\varepsilon} \left(-\frac{R^3}{k}\delta v+ \frac{R^3}{k}\,\delta(A-\alpha)\,\delta a\right)\, e^{-\int_{0}^{\hat{z}} C(\tilde{z}) d\tilde{z}} d\hat{z} \; .
\end{align*}
The change of variables $\hat{a}=A(\hat{z})$ implies 
\begin{multline*}
e^{-\int_{0}^{\varepsilon} C(\tilde{z})
	d\tilde{z}}\,\delta\rho(\varepsilon)-e^{-\int_{0}^{-\varepsilon} C(\tilde{z})
	d\tilde{z}}\,\delta\rho(-\varepsilon)=\\=\int_{A(-\varepsilon)}^{A(\varepsilon)} \left(-\frac{R^3}{k}\delta v+ \frac{R^3}{k}\,\delta(\hat{a}-\alpha)\,\delta a\right)\,\frac{1}{A'(A^{-1}(\hat{a}))} e^{-\int_{0}^{A^{-1}(\hat{a})} C(\tilde{z}) d\tilde{z}} d\hat{a} \; .
\end{multline*}

Note that $A$ is monotonically increasing with $A(0)=\alpha$. Therefore $\int_{A(-\varepsilon)}^{A(\varepsilon)} \delta(\hat{a}-\alpha)d \hat{a}=1$ and we obtain in the limit as $\varepsilon \to 0$ that
\begin{equation*}
\lim_{\varepsilon \to 0^{+}}\delta\rho(\varepsilon)-\lim_{\varepsilon \to 0^{-}}\delta\rho(\varepsilon)=\frac{R^3}{k}\,\delta a \frac{1}{A'\,(A^{-1}(\alpha))}
=\frac{R^3}{k}\,\delta a \frac{1}{A'\,(0))}  \; .
\end{equation*}
	
\section{Stability Analysis of polarisation wave S2}
\label{sec_stabanalysisS2}
Here we investigate the linear stability of the travelling wave solution S2. Linearisation of the system \eqref{equ_main} in the moving coordinate frame $(t, z=x-s t)$ at the travelling wave profile S2 lead to \eqref{equ_linearised} where $s=s_2$ as well as $R=R_2$ and $A=A_2$. The associated eigenvalue problem $(\mathcal{L}-\lambda) (\delta \rho, \delta a)^T=0$ is given by \eqref{final_A}, \eqref{equ_V}.

At the far left and right ends the travelling wave profiles $R_2$ and $A_2$ converge to
$(R_{-\infty},A_{-\infty})=(\frac{s_2}{s_2-1},1)$ and $(R_\infty,A_\infty)=(1,0)$ which corresponds to the limit values of R1 and A1 ''swapped around''.
This implies that the asymptotic matrices $B^\pm(\lambda)=\lim_{z \to \pm \infty} A(z,\lambda)$ for S2 are the same as for S1 given in \eqref{eq6_S2}, however with $+$ and $-$ inverted.
\begin{equation}\label{eq6_S2}
B^-(\lambda)=\begin{pmatrix} 0&&\frac{-s^3}{\kappa\,(s-1)^3}&&0\\
(\frac{1}{s}-1)\,\lambda&& \frac{s^2}{\kappa\,(1-s)}&& 0\\
0&& \frac{1}{1-s}&& \frac{\lambda+1}{s-1}
\end{pmatrix}\quad \text{and} \quad
B^+(\lambda)=\begin{pmatrix} 0&&\frac{-1}{\kappa}&&0\\
-\lambda&& -\frac{s}{\kappa}&& 0\\
0&& -\frac{1}{s}&& \frac{\lambda+1}{s}\end{pmatrix},
\end{equation}
where $s=s_2$. 
As a consequence the essential spectra as well as the absolute spectra of S2 and S1 coincide.

Note that the spatial eigenvalues of $B^-$ are given by those of $A^+$ given in \eqref{equ_muplus} and those of $B^+$ are given by those of $A^-$ listed in \eqref{equ_muminus}, again where $s=s_2(=-s_1)$.

To compute the Evans function, we start by identifying the stable and unstable eigenvectors of the asymptotic matrices $B^{\pm}(\lambda)$ respectively in an attempt to construct an integrable function which plays the role of an eigenvector for the linearised system. 
The spatial eigenvalues \eqref{equ_muminus} and \eqref{equ_muplus} are also spatial eigenvalues of $B^\pm$. Since here $s=s_2>1$ to avoid violations of the impenetrability constraint (see section~\ref{sec_unphysical}), the unstable eigenvalue of $B^{-}$ and the stable eigenvalues of $B^{+}$  are given by
\begin{equation}\label{eq12_2}
\mu_1^{-}(\lambda)= \frac{\lambda+1}{s-1}\;, \quad 
\mu_2^{-}(\lambda)= \frac{s^2-s \sqrt{s^2+4\kappa \lambda}}{2\kappa(1-s)}  
\quad \text{and} \quad 
\mu_3^{+}(\lambda)=\frac{-s-\sqrt{s^2+4\kappa \lambda}}{2\kappa} \; ,
\end{equation}
which satisfy $\mu_3^{+}(\lambda)<0$ and $\mu_1^{-}(\lambda)\,,\; \mu_2^{-}(\lambda)>0$ for ${\rm Re} \, \lambda \gg 1$.

We introduce the eigenvectors $X_0^+(\lambda)$ associated to the eigenvalue $\mu_2^{+}(\lambda)$ of $A^+$ as well as $X_0^-(\lambda)$ and $Y_0^-(\lambda)$ being eigenvectors of $A^-$ associated to $\mu_1^{-}(\lambda)$ and $\mu_3^{-}(\lambda)$.

To construct an eigenfunction of the linearised system we solve \eqref{final_A} using these vectors as initial, respectively terminal condition until $z=0$. Then 
the Evans function is defined as the Wronskian
\begin{equation}\label{equ_evans2}
D(\lambda)={\rm det}\; [X^{+}(z=0,\lambda),X^{-}(z=0,\lambda), Y^{-}(z=0,\lambda)] \; ,
\end{equation}
which vanishes if the stable and unstable eigenvectors propagated to $z=0$ are linearly dependent and can be combined into a smooth eigenfunction.
We start by computing $X^-(z=0, \lambda)$ and $Y^-(z=0, \lambda)$ which can be in terms of a closed-form expression. Only for the computation of $X^+(z=0, \lambda)$ we will resort to numerical results.

For $z<0$ it holds that $R  \equiv s_2/(s_2-1)$ and the matrix \eqref{final_A} is given by \eqref{eqevan1}.
Again, due to the zeros in the third column the equations for $\delta \rho$ and $\delta v$ are not coupled to $\delta a$. They satisfy \eqref{eqevan2}
which is a system of two linear, constant-coefficient equations. It admits one fundamental solution with positive eigenvalue (we omit the stable fundamental solution) given by 
\begin{equation}\label{eqevan3_S2}
\begin{pmatrix}
\delta\rho\\
\delta v
\end{pmatrix} = C_{2}\,\exp(\mu_{2}^- (\lambda) \,z ) \begin{pmatrix}
-\,\frac{s^2\,(s+\sqrt{s^2+4\,\kappa \,\lambda})}{(s-1)^2\,\kappa}\\
2\,\lambda
\end{pmatrix}
\end{equation}

With this information we can rewrite the equation for $\delta a$ which is contained in
the third row of \eqref{eqevan1},
\begin{align*}
\delta a '&=\frac{A'-1}{s-1}\,\delta v+\frac{\lambda+1}{s-1}\, \delta a\\
&=\frac{A'-1}{s-1}\,2\,\lambda\,C_{2}\,e^{\mu_{2}^-\,z} +\mu_1^-\, \delta a \; .
\end{align*}
The general solution is given by 
\begin{displaymath}
\delta a = C_1 e^{\mu_1^- z}+ \frac{2 \lambda}{s-1} C_2 \int_0^z e^{(\mu_2^--\mu_1^-) z}  \left(\frac{\alpha-1}{s-1}e^{\frac{z}{s-1}}-1 \right) \, d\tilde z \; .
\end{displaymath} 
This implies that 
$$
\lim_{z\to 0^{-}} \begin{pmatrix}
\delta\rho\\
\delta v\\
\delta a
\end{pmatrix}=\begin{pmatrix}
C_2\,\left(\,-\,\frac{s^2\,(s+\sqrt{s^2+4\,\kappa\,\lambda})}{(s-1)^2\,\kappa}\right)\\C_2\,
2\,\lambda\\
C_1
\end{pmatrix} \; .
$$ 

While most components of $A(z,\lambda)$ in \eqref{final_A} are functions, the derivative of the active speed of migration \eqref{equ_M} is a $\delta$-distribution, $M'(a)=\delta(a-\alpha)$. The change of variables between $z$ and the monotonically decreasing function $a=A(z)=A_2(z)$ shows that when integrating with respect to $z$ the following expression is a $\delta$-distribution centred at $z=0$, $-M'(A(z)) A'(z)=\delta (z)$. Therefore $\delta \rho$ and $\delta v$ in the solution of \eqref{final_A} undergo a jump at $z=0$ which involves the factor $-1/A'(0)$. The right limit as $z\rightarrow 0$ of the solution vector is given by 
\begin{align*}
\lim_{z\to 0^{+}} \begin{pmatrix}
\delta\rho\\
\delta v\\
\delta a
\end{pmatrix}=&\begin{pmatrix}
C_2\,\left(\,-\frac{s^2\,(s+\sqrt{s^2+4\,\kappa \,\lambda})}{(s-1)^2\, \kappa }\right)\\C_2\,
2\,\lambda \\
C_1
\end{pmatrix}+
\frac{(-1)C_1}{A'(0)} \frac{1}{\kappa} \begin{pmatrix}
R^3(0)\\ R(0)\,s\\
0
\end{pmatrix}
\\
=& C_2\,\begin{pmatrix}
-\,\frac{s^2\,(s+\sqrt{s^2+4\,\kappa \,\lambda})}{(s-1)^2\,\kappa}\\
2\,\lambda\\
0
\end{pmatrix}+\,C_1\,\begin{pmatrix}
-\frac{s^3}{(s-1)^2\,(\alpha-1)\,\kappa}\\
-\frac{s^2}{(\alpha-1)\,\kappa}\\ 1
\end{pmatrix}.
\end{align*} 
where we take the values of $R(0)$ and $A'(0)$ from the travelling wave profile S2.

We obtain the closed-form solutions $X^{-}(0,\lambda)$ and $Y^{-}(0,\lambda)$ corresponding to the stable eigenvalues $\mu_1^{-}(\lambda)$ and $\mu_2^{-}(\lambda)$, 
\begin{align}
\label{equ_XY_2}
X^{-}(0,\lambda) =\begin{pmatrix}
\frac{s^3}{(s-1)^2\,(\alpha-1)\,k}\\
\frac{s^2}{(\alpha-1)\,k}\\ 1
\end{pmatrix} 
\quad \text{and}\quad
Y^{-}(0,\lambda)=\begin{pmatrix}
-\,\frac{s^2\,(s+\sqrt{s^2+4\,k\,\lambda})}{(s-1)^2\,k}\\
2\,\lambda\\
0
\end{pmatrix}.
\end{align}
Finally, in order to evaluate the Evans function \eqref{equ_evans} for a given $\lambda \in \mathbb{C}$, we compute $X^+(0,\lambda)$ propagating the eigenvector associated to the spatial eigenvalue $\mu_3^+$ given by 
$$
v_3=\begin{pmatrix}
\frac{-s^2+\sqrt{s^4+4s^2\kappa\lambda}}{2\lambda s \kappa}\\
1\\ \frac{2\kappa}{s^2+2\kappa(1+\lambda)+\sqrt{s^4+4s^2\kappa\lambda}}
\end{pmatrix} \; \text{for} \; \lambda \neq 0 
\quad \text{and} \quad 
v_3=
\begin{pmatrix}
\frac{1}{s}\\
1\\ \frac{\kappa}{s^2+\kappa} 
\end{pmatrix}  \; \text{for} \; \lambda=0
$$
from an arbitrary small value of $z$ (we choose $z=20$) until $z=0$ solving the system \eqref{final_A} numerically.

\begin{figure}[H]
	\centering
	\includegraphics[width=.75\textwidth]{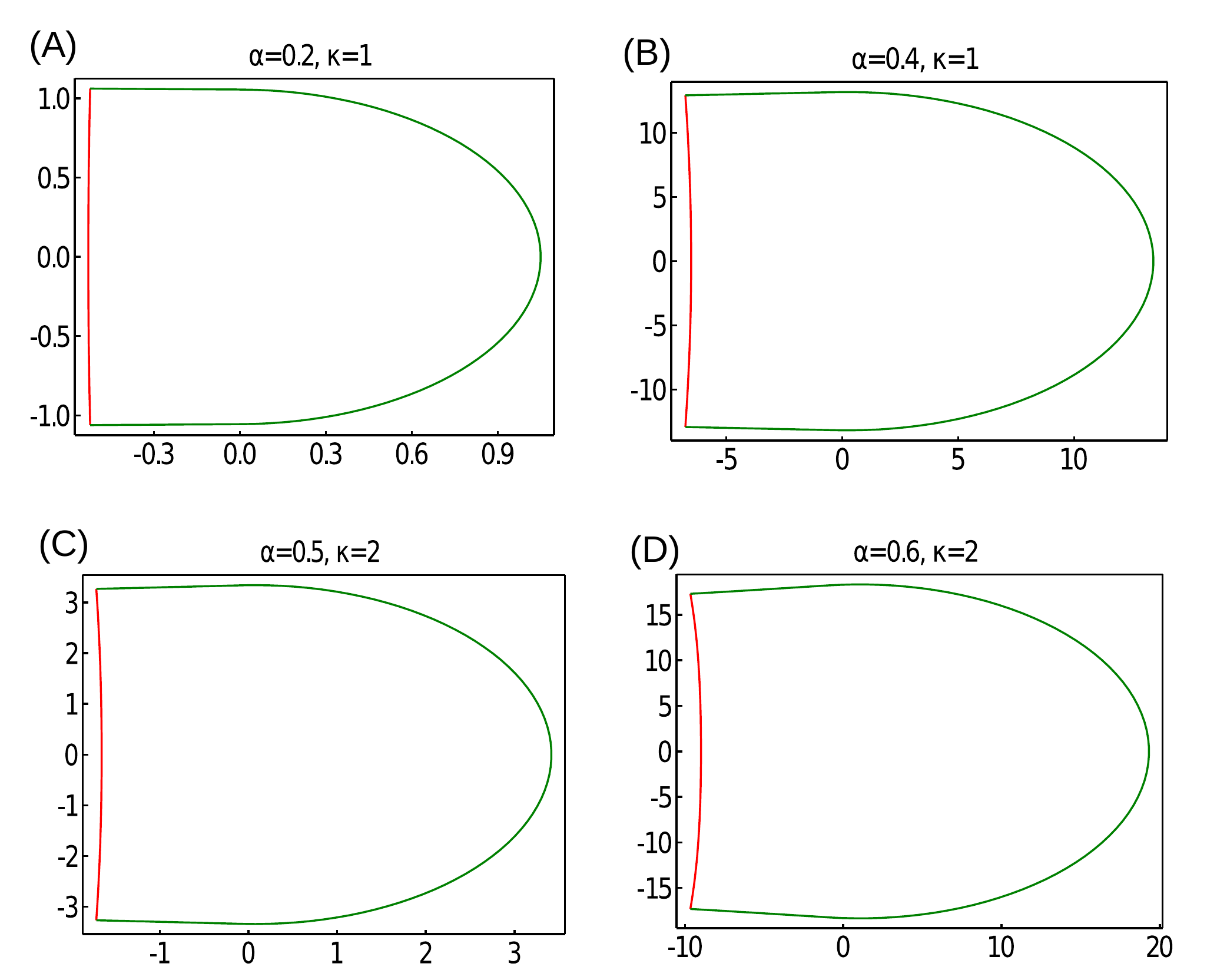}
	\caption{$D(C_1)$ for $\alpha=0.2,\,0.4,\,0.5 \quad \text{and}\quad 0.6$.}
	\label{contour1}
\end{figure}
	
\begin{figure}[H]
	\centering
	\includegraphics[width=.75\textwidth]{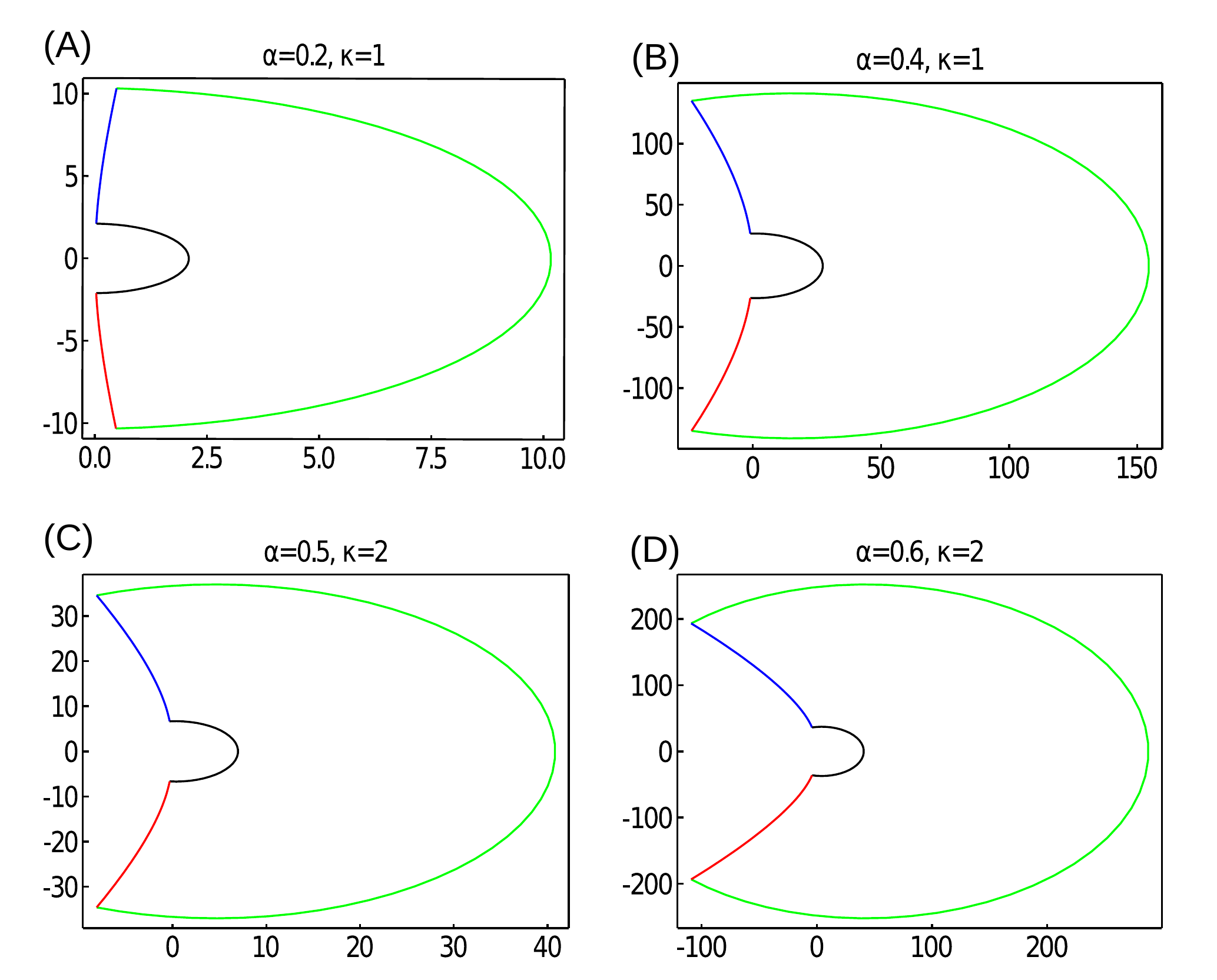}
	\caption{$D(C_2)$ for $\alpha=0.2,\,0.4,\,0.5 \quad \text{and}\quad 0.6$.}
	\label{contour2}
\end{figure} 
Finally we compute the image of both contours $C_1$ and $C_2$ (Fig.\ref{contour}) under the Evans function defined by \eqref{equ_evans2}. For various parameter values  $\alpha$ and $\kappa$ we find that for $C_1$ the winding number around the origin is $1$ (Fig.~\ref{contour1}) and for the contour $C_2$ which does not enclose the origin, the winding number is 0 (Fig.~\ref{contour2}). As for S1 this indicates that the  point spectrum on the right of the complex plane only consists of $\lambda=0$ which is expected for a travelling wave solution. This suggests that the travelling wave solutions S2 and -- through to transformation $T_2$ -- S4 are linearly stable. 

\section{Unphysical travelling wave solutions}
\label{sec_unphysical}

\begin{figure}[H]
	\centering
	\includegraphics[width=.5\textheight]{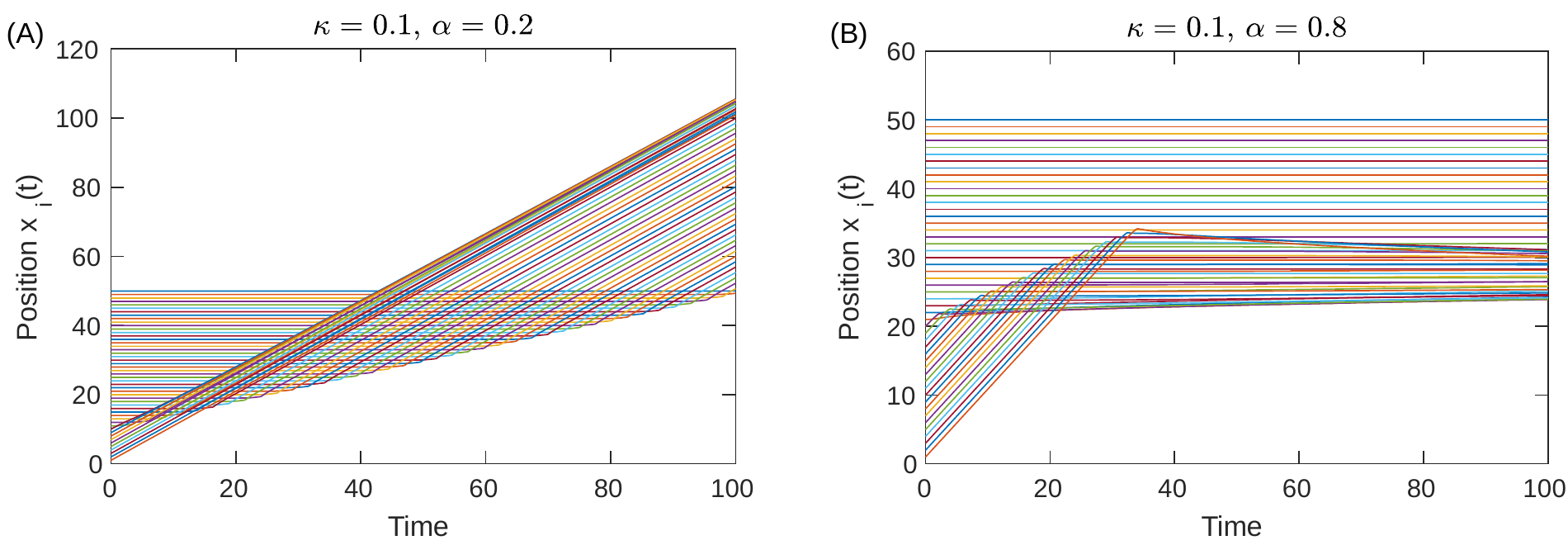}
	\caption{Un-physical travelling wave solutions in violation of the cells' impenetrablity. (A) shows a simulation of the polarisation wave S2 for a set of parameters for which $s_2<1$. (B) shows a simulation of the polarisation wave S4 for a set of parameters for which $s_4>0$.}
	\label{unphysical}
\end{figure}

Note that the travelling wave solution S2 can only be realised if the model parameters $\kappa$ and $\alpha$ are such that $s_2>1$. If that is not the case the mathematical solution is not physical and violates the impenetrability of single cells as illustrated in Fig.~\ref{unphysical}.

\bibliographystyle{IEEEtran}
\bibliography{sample}
\end{document}